\documentclass[11pt]{article}

\usepackage{graphicx}
\usepackage{amsmath}
\usepackage{amsfonts}
\usepackage{fancyhdr} 
\usepackage{mathrsfs}

\usepackage{wrapfig}

\usepackage{bbold}

\newtheorem{lem}{Lemma}
\newtheorem{thm}{Theorem}

\newtheorem{prop}[lem]{Proposition}
\newtheorem{cor}[lem]{Corollary}
\newtheorem{defn}[lem]{Definition}

\newenvironment{thmbis}[1]
  {\renewcommand{\thethm}{\ref{#1}$'$}%
   \addtocounter{thm}{-1}%
   \begin{thm}}
  {\end{thm}}

\numberwithin{equation}{section}
\numberwithin{thm}{section}
\numberwithin{lem}{section}
\numberwithin{alg}{section}

\DeclareMathOperator{\sinc}{sinc}

\newcommand{\eee}{\mathcal{E}}
\newcommand{\fff}{\mathcal{F}}
\newcommand{\sss}{\mathscr{S}}

\newcommand{\real}{\mathbb{R}}
\newcommand{\integer}{\mathbb{Z}}

\newcommand{\ddd}{\mathbb{D}}

\newcommand{\complex}{\mathbb{C}}

\newenvironment{pf}{\noindent {\em Proof}.\ \ }{\hspace*{\fill}\rule{.5ex}{1.4ex}\,}

\newcommand{\lll}{\mathcal L}
\newcommand{\one}{\mathbb{1}}
\newcommand{\hhh}{\mathscr{H}}
\newcommand{\ggg}{\mathscr{G}}

\newcommand{\ccc}{\mathcal{C}}

\newcommand{\ttt}{\mathbb{T}}

\newcommand{\dist}{\mathcal{D}^\prime}
\newcommand{\kkk}{\mathscr{K}}
\newcommand{\mlap}{\widetilde{\Delta}}

\DeclareMathOperator{\spn}{span}

\title{\vspace*{-1in}A smooth global model for scattering in layered media}
\author{Peter C.~Gibson\footnote{Dept.~of Mathematics \& Statistics, York University, 4700 Keele St., Toronto, Ontario, Canada, M3J~1P3, $\mathtt{pcgibson@yorku.ca}$} }
\date{December 18, 2014}
\makeatletter
\let\newtitle\@title
\let\newauthor\@author
\let\newdate\@date
\makeatother

\begin{document}
\maketitle
\begin{abstract}
Layered media have been studied extensively both for their importance in imaging technologies and as an example of a hyperbolic PDE with discontinuous coefficients.   From the perspective of acoustic imaging, the time limited impulse response at the boundary, or boundary Green's function, represents measured data, and the objective is to determine coefficients, which encode physical parameters, from the data.   The present paper resolves two fundamental open problems for layered media: (1) how to compute the time limited Green's function in the presence of discontinuous coefficients; and (2) to determine precisely how data depends on coefficients.  We show that there exists a single system of equations in $3n$-dimensional space that governs the parameterized family of all $n$-layered media simultaneously. The alternate system has smooth coefficients, can be solved directly by separation of variables, and recovers the impulse response at the boundary for the original equations. The analysis brings to light an exotic laplacian---hybrid between the euclidean and hyperbolic laplacians---that plays a central role in the scattering process.  Its eigenfunctions comprise a new family of orthogonal polynomials on the disk.  These serve as building blocks for a universal wavefield in terms of which the dependence of data on coefficients has a simple description: reflection data is obtained by sampling a translate of the wavefield on the integer lattice and then pushing forward by a linear functional, where the translate and pushforward correspond to reflectivity and layer depth vectors, respectively.   
\end{abstract}

\tableofcontents

\section{Introduction}

\pagestyle{fancyplain}

The classical model for scattering in layered media, detailed below in \S\ref{sec-standard-model}, is important both in applications to imaging technologies and as a theoretical case study.  Piecewise constant layered media serves as a basic model in seismic imaging \cite{BlCoSt:2001},\cite{Yi:2001}, historically the main driver of research into the subject, as well as in acoustic and electromagnetic imaging \cite{KuFe:2009}, and the design of optical coatings \cite{FuTi:1992}, for example.  From the theoretical point of view, the governing equations for piecewise constant layered media comprise a basic example of a wave equation with discontinuous coefficients, a feature that is incompatible with a wide swath of established theory (e.g., \cite{Bu:1980}, \cite{Sy:1983}, \cite{Ra:2003}, \cite{Ra:2008}).  In imaging applications, the impulse response at the boundary for these equations, or boundary Green's function, corresponds to measured data, which can be recorded only for a finite length of time.  A central theoretical difficulty stems from the combination of discontinuous coefficients with finite time duration; indeed, despite the extensive literature on layered media, an explicit closed-form formula for the time limited Green's function was discovered only recently \cite{Gi:SIAP2014}.  Consequently, the way measured data depends on physical parameters, encoded as coefficients of the governing equations, has remained obscure, allowing only a partial understanding of the associated inverse problem.   

The purpose of the present paper is to present a novel perspective on piecewise constant layered media that resolves the central theoretical difficulty and reveals an underlying mathematical structure not previously known.  We establish several basic new results, summarized as follows. 
\renewcommand{\theenumi}{\Roman{enumi}}
\begin{enumerate}
\item There is a single global system of first and second order PDEs that governs the impulse responses at the boundary for all $n$-layered media collectively.   
\item The global system has smooth coefficients and can be solved directly by classical separation of variables to yield explicit formulas for the impulse responses in a form suited to time-limited data. 
\item There is a universal amplitude wavefield $\psi$ on $\complex^n$ in terms of which the dependence of reflection data on material parameters has a simple interpretation:  data is realized by sampling a translate of $\psi$ on the integer lattice and then pushing this forward by a linear functional (onto a one-dimensional timeline), where the translate and pushforward correspond respectively to reflectivity and layer depth vectors.  
\item Eigenfunctions of the two-dimensional modified laplacian 
$
\frac{1-x^2-y^2}{4}\Delta
$
comprise a previously unknown family of orthogonal polynomials on the disk that serve as building blocks for the universal amplitude wavefield.  
\end{enumerate}
\renewcommand{\theenumi}{\arabic{enumi}}
These results show for the first time that the nonlinear dependence of data on material parameters is governed by a system of linear PDEs.  Our analysis brings to light a particular riemannian metric on the disk, corresponding to the above modified laplacian, that is central to the paper's results; the smooth global model and its solution are formulated in terms of the metric's Laplace-Beltrami operator and its eigenfunctions, respectively.   From the perspective of general theory, the results give an indication of what the correspondence between data and physical parameters might look like in a more general setting; for instance, although the correspondence is nonlinear, it makes sense to ask whether it is governed by a system of linear equations.  

\subsection{Related literature\label{sec-related}}
Mathematical literature on acoustic imaging of layered media has accrued steadily over the last half century, driven by applications to seismic exploration, ultrasound and other technologies.  A standard approach that dates to the early 1950s, \cite{Br:1951}, is to patch together solutions within each of the layers, imposing continuity conditions at layer boundaries; see \cite[Ch.~3]{FoGaPaSo:2007}, \cite[Ch.~6]{BrGo:1990}, \cite{BuBu:1983} and the many references therein.  Reflection and transmission at layer boundaries engender a cascade of progressively more complicated scattering series which are in practice approximated as in \cite{In:2009}.  Such approximations are avoided altogether in the present paper; the global model eliminates the need to patch together solutions, instead making exact results accessible by standard methods.    

It is well known that data depends nonlinearly on physical parameters for a variety of inverse problems in mathematical physics, including gravimetry, conductivity and tomography \cite[Ch.~1]{Is:1998}.  Because the precise nature of the nonlinearity is in some cases unknown, the linearized correspondence is considered instead.  This leads to a Newton type approach to the inverse problem, a rigorous mathematical framework for which has been developed in \cite{BlStSy:2013}, for a broad class of hyperbolic equations.   Newton type methods have the drawback that an a priori initial estimate for unknown coefficients is required.  By giving a precise description of the nonlinear dependence of data on coefficients, the present paper shows that there is an alternative to linearization in the case of the classical model for layered media.  Thus the need for a priori information can be avoided, which is highly desirable from the perspective of imaging.  
%It remains to be seen whether this can be extended to a wider class of models.  

\section{The usual model for layered media\label{sec-standard-model}}

\begin{figure}[h]
\hspace*{\fill}
\parbox{5.8in}{
\fbox{
\includegraphics[clip,trim=10in 0in 15in 0in, width=5.8in]{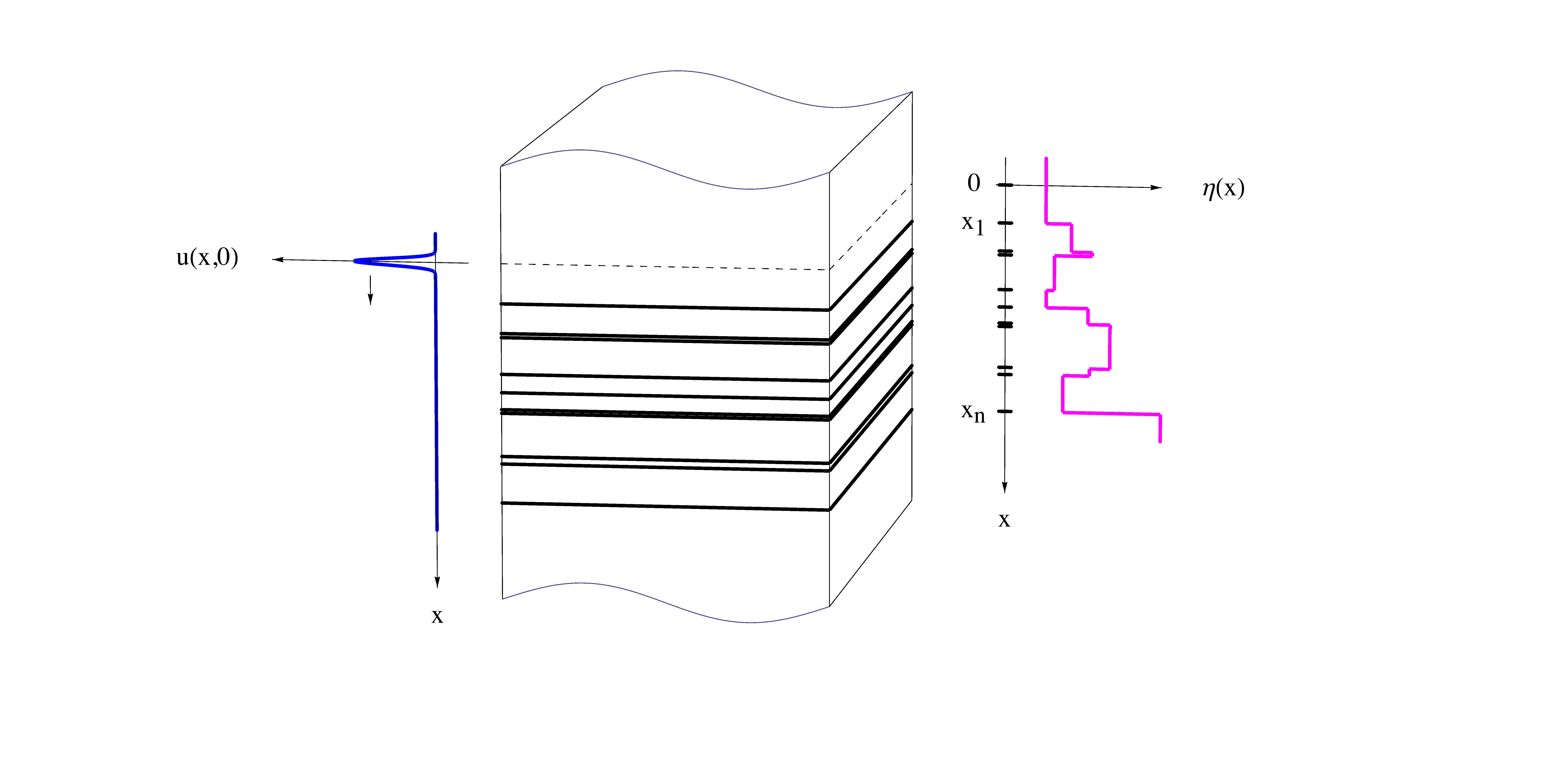}
%left,bottom,right,top
}
\caption{A layered medium.  An impulsive source wave (blue) travels from the source depth $x=0$ (dashed line) toward the layers.  In acoustic imaging the objective is to recover the impedance profile (magenta) from time-limited measurements at the source depth of the scattered wave.}\label{fig-experiment}
}\hspace*{\fill}
\end{figure}

Consider a three-dimensional medium consisting of $n$ homogeneous horizontal layers sandwiched between two half spaces.  Suppose that a horizontal plane wave impulse is transmitted toward the layers from a source depth above the layers, and that the resulting scattered wavefield is recorded at the source depth over some finite time interval $[0,T]$.   
This scenario is modeled by a hyperbolic system
\begin{align}\label{model}
\begin{split}
\eta\frac{\partial u}{\partial t}+\frac{\partial p}{\partial x}&=\eta^{1/2}(0)\delta(t)\delta(x)\\
\frac{\partial p}{\partial t}+\eta\frac{\partial u}{\partial x}&=0
\end{split}
%\begin{split}
%u(0,x)&=\delta(x)\\
% p(0,x)&=\eta(0)\;\delta(x)
%\end{split}
\end{align}
%$\eta\in\fff_n^\real$
where $u(t,x)$ and $p(t,x)$ denote particle velocity and pressure, respectively, as functions of time $t$ and renormalized depth $x$, and where $\eta(x)$ denotes acoustic impedance as a function of renormalized depth.  See Figure~\ref{fig-experiment}.  

Assuming that the system is initially quiescent, meaning $u(t,x)=p(t,x)=0$ for $t<0$, let $G_\eta(t)$ denote the velocity part of the unique distributional solution to the system (\ref{model}), restricted to source depth $x=0$:  
\begin{equation}\label{response}
G_\eta(t)=u(t,0).
\end{equation}
The scattered wavefield recorded at source depth $x=0$ is 
\begin{equation}\label{data}
\chi_{[0,T]}G_\eta,
\end{equation}
which we shall refer to as the measured data.   The forward problem is to determine the data $\chi_{[0,T]}G_\eta$ given a particular impedance profile $\eta$, and the inverse problem is to determine $\eta$ given $\chi_{[0,T]}G_\eta$.   The impulse response is easily seen to have the structure of a delta train
\[
G_\eta(t)=\sum_{j=1}^\infty a_j\delta(t-t_j),
\]
of which the measured data is a partial sum
\begin{equation}\label{measured}
\chi_{[0,T]}G_\eta(t)=\sum_{j=1}^Na_j\delta(t-t_j)
\end{equation}
for some $N$ that depends both on $\eta$ and $T$.  The forward problem is thus to express the sequences
\begin{equation}\label{amplitudes-arrival-times}
(a_1,\ldots,a_N)\quad\mbox{ and }\quad(t_1,\ldots,t_N)
\end{equation}
in terms of $\eta$ and $T$, and the inverse problem to recover $\eta$ from (\ref{amplitudes-arrival-times}).

The model (\ref{model}) is classical and has been studied extensively; see \cite[Chapter~3]{FoGaPaSo:2007} for a derivation from physical principles and an up-to-date summary of what is known.  Connections to various other models appear, for example, in \cite{Bu:1980}, as do details concerning renormalization of depth.

\subsection{The complexified model\label{sec-complexified-model}}

From a physical perspective acoustic impedance is a nonnegative quantity.  We complexify the classical model by allowing impedance to take complex values with nonnegative real part, a modification which turns out to be key.  The complexified impedance transforms into a sequence of reflection coefficients taking values in the closed disk, and we show in later sections that there is an associated riemannian structure on the disk that illuminates both the forward and inverse problems for the original model. This structure is not visible without complexification.    

We fix notation to make precise the family of media under consideration, as follows.  Let $H$ denote the Heaviside function. 
%\begin{equation}\label{heaviside}
%H:\real\rightarrow\{0,1\},\quad\quad H(x)=\left\{\begin{array}{cc}0&\mbox{ if }x\leq 0\\1&\mbox{ if }x>0\end{array}\right..
%\end{equation}
The right half plane, open unit disk and unit circle are respectively denoted
\[
\complex_+=\{z\in\complex\,|\,\Re z>0\},\quad\ddd=\{z\in\complex\,|\;|z|<1\},\quad
\ttt=\{z\in\complex\,|\;|z|=1\}.
\]
Generally speaking $n\geq 2$ denotes the number of layers in layered media under consideration.  Define $\fff_n$ to be the collection of all functions
\[
\eta:\real\rightarrow\complex_+
\]
of the form
\begin{equation}\label{eta}
\begin{split}
&\eta(x)=C_0+\sum_{j=1}^nC_jH(x-X_j)\\
\mbox{ where }\quad
&0<X_1<X_2<\cdots<X_n\quad\mbox{ and }\quad\sum_{j=0}^nC_j=1.
\end{split}
\end{equation}
Note that each $\eta\in\fff_n$ is defined to have values in the right half plane $\complex_+$, so that 
\[
(C_0,C_0+C_1,\ldots,C_0+C_1+\cdots+C_{n-1})\in\complex_+^n.
\]
We write $\fff_n^\real\subset\fff_n$ to denote the subset consisting of real-valued functions $\eta:\real\rightarrow\real_+$, which correspond to physically realizable systems.\footnote{The normalization $\sum C_j=1$ selects a single impedance profile to represent the equivalence class of media having a common impulse response $G_\eta$.}

Viewed globally, the model (\ref{model}) consists of an uncountably infinite family of coupled equations parameterized by complexified impedance $\eta\in\fff_n$.   A principal contribution of the present paper is to show this infinite family of coupled equations can be replaced by a single system of $2n$ equations that have smooth coefficients---a smooth global model.  
The smooth model makes it possible to derive an explicit formula for the measured data $\chi_{[0,T]}G_\eta$ directly by separation of variables.  By contrast, the usual model (\ref{model}) is not amenable to such a straightforward analysis.  Indeed, the discontinuous coefficients obstruct not only separation of variables, but also any other techniques that require at least $H^1$ regularity (e.g., \cite{Bu:1980}, \cite{Sy:1983}, \cite{Ra:2003}, \cite{Ra:2008}).   

On the other hand, there is an elegant and well-known expression for the Fourier transform of $G_\eta$---without the cutoff function $\chi_{[0,T]}$---that can be derived easily in the traditional setting of real-valued impedance \cite[\S3.5]{FoGaPaSo:2007}.  This expression is a backward recurrence formula, and it too naturally lends itself to complexification.  In the next sections we use a complexified version of the backward recurrence to describe the essential difficulty inherent in computing $\chi_{[0,T]}G_\eta$.

\subsection{The complexified backward recurrence formula\label{sec-backward-recurrence}}

The present paper uses the version of the Fourier transform consistent with the formula
\begin{equation}\label{Fourier}
\hat{f}(\sigma)=\int_{\real}f(t)e^{i\sigma t}\,dt.  
\end{equation}
To formulate the backward recurrence we convert the impedance $\eta$ into reflectivity and layer depth, for which purpose we pass between $\complex_+^n$ and the open polydisk by way of the bijective mapping 
\[
\Phi:\complex_+^n\rightarrow\ddd^n
\]
defined by the formula
\begin{equation}\label{Phi}
\Phi_j(\zeta)=\left\{\begin{array}{cc}
(\zeta_j-\zeta_{j+1})/(\zeta_j+\overline{\zeta_{j+1}})&\mbox{ if }1\leq j\leq n-1\\[9pt]
(\zeta_n-1)/(\zeta_n+1)&\mbox{ if }j=n
\end{array}\right. .
\end{equation}
For a given impedance profile
\[
\eta(x)=C_0+\sum_{j=1}^nC_jH(x-X_j),
\]
set
\begin{equation}\label{reflectivities}
(w_1,\ldots,w_n)=\Phi(C_0,C_0+C_1,\ldots,C_0+C_1+\cdots+C_{n-1})
\end{equation}
and 
\begin{equation}\label{travel-times}
(\tau_1,\ldots,\tau_n)=(X_1,X_2-X_1,\ldots,X_n-X_{n-1}).
\end{equation}
In the case where $\eta\in\fff_n^\real$, entries of the resulting vector $w=(w_1,\ldots,w_n)$ are reflection coefficients at layer boundaries, and $\tau=(\tau_1,\ldots,\tau_n)$ is the vector of layer depths in units of two-way travel time; the formula (\ref{Phi}) expresses the standard transformation from impedance to reflectivity.    

For each pair $(w,z)\in\overline{\ddd}^n\times\ttt^n$, and each $v\in\overline{\ddd}$, write 
\begin{equation}\label{Psi}
\Psi_{z_j}^{w_j}(v)=z_j\frac{v+\overline{w}_j}{1+w_jv}\quad\quad(1\leq j\leq n).
\end{equation}
In terms of this notation the complexified backward recurrence formula is 
\begin{equation}\label{backward}
\widehat{G}_\eta(\sigma)=\Psi_{e^{i\sigma\tau_1}}^{w_1}\circ\Psi_{e^{i\sigma\tau_2}}^{w_2}\circ\cdots\circ\Psi_{e^{i\sigma\tau_n}}^{w_n}(0).
\end{equation}
The real version of this (without the conjugation that appears in (\ref{Psi})) has been known since the earliest work on the subject, cf., \cite{Br:1951}, \cite{Ku:1963}.

Computation of $\widehat{G}_\eta(\sigma)$ by means of the formula (\ref{backward}) is fast and exact. However the situation changes dramatically with the introduction of a cutoff function $\chi_{[0,T]}(t)$.   
The Fourier transform of $\chi_{[0,T]}G_\eta$ is 
\begin{equation}\label{truncated-transform}
K_T\ast\widehat{G}_\eta,\quad\mbox{ where }\quad K_T(\sigma)=\textstyle\frac{T}{2\pi}e^{i\sigma T/2}\sinc(\sigma T/2).  
\end{equation}
Unlike the pure backward recurrence, this cannot be computed quickly or exactly, since $K_T$ is supported on the whole real line, and the convolution requires knowing $\widehat{G}_\eta$ everywhere, not just at a single point.  Thus from a computational perspective the backward recurrence does not yield an exact solution to the forward problem.

\subsection{Essential questions\label{sec-essential-questions}}

The Fourier transform of the measured data represented in the form (\ref{measured}) is 
\begin{equation}\label{Fourier-measured}
K_T\ast\widehat{G}_\eta(\sigma)=\sum_{j=1}^Na_je^{it_j\sigma},
\end{equation}
which is a truncation of the almost-periodic expansion of $\widehat{G}_\eta$ itself,
\begin{equation}\label{almost-periodic}
\widehat{G}_\eta(\sigma)=\sum_{j=1}^\infty a_je^{it_j\sigma}.
\end{equation}
Truncation being a trivial matter, the essential problem is to compute explicitly the almost-periodic expansion of the Fourier transform of $G_\eta$, i.e., of the backward recurrence formula (\ref{backward}), 
\begin{equation}\label{expansion}
\Psi_{e^{i\sigma\tau_1}}^{w_1}\circ\cdots\circ\Psi_{e^{i\sigma\tau_n}}^{w_n}(0)=\sum_{j=1}^\infty a_je^{it_j\sigma}.
\end{equation}
This is the heart of the matter from the analytic point of view.  We are searching for an analogue of the various well-known expansions in terms of Bessel functions, Hankel functions or orthogonal polynomials that arise in other types of scattering \cite[Ch.~6]{Is:1998}, \cite{CoKr:2013}.  Such an expansion will enable an explicit representation of the measured data in the time domain, and facilitate a direct analysis of a second essential question of how precisely measured data depends on coefficients.   Almost periodicity of the right-hand side of (\ref{expansion}) is a technical complication that means that the analogy cannot be exact; some extra structure must come into play.  

It turns out that there do indeed exist special functions appropriate to the right-hand side of (\ref{expansion}).  (Their existence was first hinted at in \cite{Gi:SIAP2014}, where an explicit formula for $\chi_{[0,T]}G_\eta$, for $\eta\in\fff^\real_n$, was computed using combinatorial means.) They comprise a new class of orthogonal polynomials on the disk (cf.~\cite{DuXu:2001} and \cite{Xu:2015}), and are related to a seldom-used riemannian metric that is intermediate between the euclidean and hyperbolic metrics.  We introduce the Laplace-Beltrami operator for this metric below, as a prelude to setting up a smooth global model.

\section{Statement of the main results\label{sec-summary}}

\subsection{Key ingredients: layer collapse and a hybrid laplacian\label{sec-key-ingredients}}

In this section we describe several key ingredients to the global model presented in the subsequent section: (1) a boundary condition realized by collapsing layers in the original model so as to wash out time dependence; (2) the Laplace-Beltrami operator for a particular metric on the unit disk; and (3) a simple pushforward construction. 

Setting  each $\tau_j=0$ in the formula (\ref{backward}) corresponds to collapsing the layers in the physical model to $n$ infinitesimally thin layers located together at the source depth $x=0$.  From the global perspective the collapsed structure comprises boundary data.  More precisely, define the \emph{layer collapse function} on the polydisk as 
\begin{equation}\label{eee-definition}
\eee:\overline{\ddd}^n\rightarrow\overline{\ddd},\qquad\eee(w)=\Psi^{w_1}_1\circ\Psi^{w_2}_1\circ\cdots\circ\Psi^{w_n}_1(0).
\end{equation}
This simplified version of the backward recurrence is much easier to analyze than the original formula, but it nevertheless encodes crucial information.  

Next we define the \emph{hybrid laplacian}, which is, roughly speaking, hybrid between the euclidean laplacian and the hyperbolic laplacian.  In the context of the disk, we use the notation $z=x+iy$ and pass between complex and euclidean coordinates without further comment.  Set 
\begin{equation}\label{hybrid-defintion}
\mlap=(1-z\bar{z})\frac{\partial^2}{\partial z\partial\bar{z}}=\frac{1-x^2-y^2}{4}\left(\frac{\partial^2}{\partial x^2}+\frac{\partial^2}{\partial y^2}\right)=\frac{1-x^2-y^2}{4}\Delta,
\end{equation}
where $\Delta$ is the usual euclidean laplacian.  (Note that the hyperbolic laplacian has the form $4\Delta/(1-x^2-y^2)^2$.)  Thus $\mlap$ is the Laplace-Beltrami operator for the riemannian metric
\begin{equation}\label{hybrid-metric}
ds^2=\frac{4}{1-x^2-y^2}\left(dx^2+dy^2\right),
\end{equation}
which degenerates at the boundary circle $\ttt$ (and which has area measure
%\begin{equation}\label{area-form}
$d\mu=4dx\,dy/(1-x^2-y^2)$.) 
%\end{equation} 
With respect to the metric (\ref{hybrid-metric}) the disk has infinite area but finite diameter.  Because of the latter property, the geodesic distance from any given point in $\ddd$ to the boundary is finite, and one can attach $\ttt$ to form the riemannian manifold with boundary
\begin{equation}\label{space-with-boundary}
\bigl(\overline{\ddd},ds^2\bigr),
\end{equation}
which we henceforth refer to as the \emph{scattering disk}.   

Given a complex $n$-tuple $w=(w_1,\ldots,w_n)\in\ddd^n$, let $\Delta_j$ denote the standard laplacian with respect to the variable $w_j$, for $1\leq j\leq n$, and let
\begin{equation}\label{modified-wj}
\widetilde{\Delta}_j=\frac{1-|w_j|^2}{4}\Delta_j  
\end{equation}
denote the corresponding hybrid laplacian.  

Lastly we describe a simple pushforward that facilitates treating the family of distributions $G_\eta$ collectively.  Convert  $\eta\in\fff_n$ into a pair $(\tau,w)$ of $n$-tuples by the transformations (\ref{travel-times}) and (\ref{reflectivities}).  Then the travel time vector $\tau$ in particular may be viewed as a linear functional on $\real^n$ (which we also denote as $\tau$) under the action of the standard scalar product
\[
\tau:\real^n\rightarrow\real, \quad x\mapsto\langle x,\tau\rangle.
\]
Viewing $\real^n$ as space-time, $\tau$ thus maps space-time onto a timeline $t=\langle x,\tau\rangle$.  Distributions $F$ on $\real^n$ push forward by $\tau$ to distributions $\tau_*F$ on $\real$ according to the formula
\begin{equation}\label{pushforward}
(\tau_\ast F,\varphi)=(F,\varphi\circ\tau),
\end{equation}
where $(\cdot,\cdot)$ denotes the pairing of a distribution with a test function,\footnote{Of course it has to be possible to give a precise interpretation to the right-hand pairing $(F,\varphi\circ\tau)$, since $\varphi\circ\tau$ is in general not a proper test function on $\real^n$.  For present purposes this requires giving a meaningful interpretation to evaluation of a Dirac delta on a smooth function that is constant on hyperplanes, which is always possible.}\! and $\varphi(t)$ is a test function on $\real$.  To exploit this idea we introduce space-time variables $x\in\real^n$ and seek a distribution $\ggg(w,x)$ that pushes forward to $G_\eta(t)$ via $\tau_\ast$, for every $\eta\in\fff_n$.  The Fourier dual notion to pushforward is restriction; in the dual domain $\widehat{\ggg}(w,\xi)$ restricts to $\widehat{G}_\eta(\sigma)$ along the line $\xi=\sigma\tau$.

\subsection{A global model\label{sec-global-model}}

Let $U:\overline{\ddd}^n\times\real^n\rightarrow\complex$ denote an unknown distribution of the general form   
\begin{equation}\label{U-form}
U(w,x)=\sum_{k\in\integer^n}c_k(w)\delta(x-k).
\end{equation}
(This condition is equivalent to periodicity of the Fourier transform with respect to $x$.)   
We write $w_j=r_je^{i\theta_j}$ for the polar form of entries of $w$ and adopt the convention that $x_{n+1}=0$.  For $1\leq j\leq n$, consider the Helmholtz-type equations
\begin{equation}\label{Helmholtz}
\widetilde{\Delta}_jU+x_jx_{j+1}U=0
\end{equation}
and the first order coupling equations
\begin{equation}\label{coupling}
\frac{\partial U}{\partial \theta_j}+i(x_j-x_{j+1})U=0
\end{equation}
subject to the trace condition
\begin{equation}\label{trace}
\int_{\real^n}U(w,x)\,dx=\eee(w)\quad\mbox{ for } w\in\partial\overline{\ddd}^n.
\end{equation}
Our main results assert that the system (\ref{U-form}), (\ref{Helmholtz}), (\ref{coupling}), and (\ref{trace}) comprises a global model for scattering in layered media that is equivalent to the family of models (\ref{model}) for $\eta\in\fff_n$.  Crucially, the global model is amenable to solution by separation of variables, without any additional machinery or ad hoc arguments, essentially because equations (\ref{Helmholtz}) and (\ref{coupling}) have smooth coefficients---unlike (\ref{model}).  And separation of variables produces precisely the type of expansion called for earlier in \S~\ref{sec-essential-questions}, expressed in terms of the following polynomials.  

\begin{defn}[scattering polynomials]\label{defn-scattering-polynomials}
For each $(p,q)\in\integer^2$ define the scattering polynomial $\varphi^{(p,q)}:\complex\rightarrow\complex$ as follows.  If $\min\{p,q\}\geq 1$, set 
\begin{equation}\label{scattering-polynomial}
\varphi^{(p,q)}(\zeta)=
\alpha_{p,q}\bigl(1-\zeta\bar{\zeta}\bigr)\displaystyle\frac{\partial^{p+q}}{\partial\zeta^p\partial\bar{\zeta}^q}\bigl(1-\zeta\bar{\zeta}\bigr)^{p+q-1},\quad\mbox{ where }\quad\alpha_{p,q}=\frac{(-1)^p}{q(p+q-1)!}.
\end{equation}
If $p\geq 0$ set $\varphi^{(p,0)}(\zeta)=\bar{\zeta}^p$; otherwise set $\varphi^{(p,q)}=0$.
\end{defn}
The scattering polynomials map the closed unit disk into itself; see \S\ref{sec-new-polynomials}.  
They have integer coefficients as polynomials in variables $\zeta,\bar{\zeta}$; expressed in terms of variables $x,y$, where $\zeta=x+iy$, they have Gaussian integer coefficients belonging to the the set $\integer\cup i\integer$.  Most importantly, they are eigenfunctions of the hybrid laplacian, as follows.  
\begin{thm}\label{thm-hybrid-laplacian}
Given Dirichlet boundary conditions on $\overline{\ddd}$, the non-zero eigenvalues of $-\mlap$ are positive integers.  For each positive integer $k$,
\[
\ker(\mlap+k)=\spn\bigl\{\varphi^{(p,q)}\,|\;pq=k\mbox{ and }p,q\in\integer_+\bigr\}.
\]
In particular, the dimension of $\ker(\mlap+k)$ is the number of divisors of $k$.  
\end{thm}

Since equation (\ref{Helmholtz}) involves the hybrid laplacian, the solution to our global model is most simply expressed in terms of scattering polynomials; in the following statement we use the convention that for $k\in\integer^n$, $k_{n+1}=0$.   
\begin{thm}\label{thm-unique}
The unique distributional solution to the system (\ref{U-form}), (\ref{Helmholtz}), (\ref{coupling}), and (\ref{trace}) is
\[
\ggg(w,x)=\sum_{k\in\{1\}\times\integer^{n-1}}\Bigl(\prod_{j=1}^n\varphi^{(k_j,k_{j+1})}(w_j)\Bigr)\delta(x-k).
\]
\end{thm}
The next result invokes the pushforward (\ref{pushforward}).   
\begin{thm}\label{thm-capture}
Convert an arbitrary $\eta\in\fff_n$ into a pair $(\tau,w)$ of $n$-tuples by the transformations (\ref{travel-times}) and (\ref{reflectivities}). Then $G_\eta(t)=\tau_\ast\ggg(w,t)$.  
\end{thm}
Explicit evaluation of the pushforward $\tau_\ast\ggg$ yields the formula 
\begin{equation}\label{explicit-pushforward}
G_\eta(t)=\sum_{k\in\{1\}\times\integer^{n-1}}\Bigl(\prod_{j=1}^n\varphi^{(k_j,k_{j+1})}(w_j)\Bigr)\delta(t-\langle k,\tau\rangle).
\end{equation}

In the Fourier domain the system (\ref{U-form}), (\ref{Helmholtz}), (\ref{coupling}), and (\ref{trace}) becomes even simpler, as follows.  Let $\widehat{U}(w,\xi)$ denote the Fourier transform of $U(w,x)$ in $x$ consistent with the formalism
\[
\widehat{U}(w,\xi)=\int_{\real^n}U(w,x)\,e^{i\langle x,\xi\rangle}\,dx.
\]
Condition (\ref{U-form}) is equivalent to periodicity of  $\widehat{U}$.  More precisely, $\widehat{U}$ is the pullback to $\overline{\ddd}^n\times\real^n$ of a function
\[
V:\overline{\ddd}^n\times\ttt^n\rightarrow\complex,
\]
with 
\begin{equation}\label{V}
\widehat{U}(w,\xi)=V(w,z),\quad z=(e^{i\xi_1},\ldots,e^{i\xi_n}).
\end{equation}
In terms of this notation, for $1\leq j\leq n$, the function $V$ satisfies the wave equations
\begin{equation}\label{wave}
-\widetilde{\Delta}_jV+\frac{\partial^2V}{\partial\xi_{j+1}\partial\xi_j}=0,
\end{equation}
and the coupling equations
\begin{equation}\label{coupling-dual}
\frac{\partial V}{\partial\theta_j}+\frac{\partial V}{\partial\xi_j}-\frac{\partial V}{\partial\xi_{j+1}}=0,
\end{equation}
subject to the partial boundary condition 
\begin{equation}\label{V-boundary}
V(w,\one)=\eee(w)\quad\mbox{ for } w\in\partial\overline{\ddd}^n,
\end{equation}
where $\one$ denotes the constant vector $\one=(1,\ldots,1)$.  
 
Keeping the notation (\ref{V}), Theorems~\ref{thm-unique} and \ref{thm-capture} have an equivalent formulation in terms of $V$, as follows. 
\begin{thmbis}{thm-unique}
The unique distribution on $\overline{\ddd}^n\times\ttt^n$ satisfying equations (\ref{wave}),  (\ref{coupling-dual}) and the partial boundary condition (\ref{V-boundary}) is
\[
\hhh(w,z)=\sum_{k\in\{1\}\times\integer^{n-1}}\Bigl(\prod_{j=1}^n\varphi^{(k_j,k_{j+1})}(w_j)\Bigr)z^k.
\]
\end{thmbis}
Note that $\hhh$ pulls back to $\widehat{\ggg}$ in the obvious way,
\[
\widehat{\ggg}(w,\xi)=\hhh(w,z),
\]
where $z=\bigl(e^{i\xi_1},\dots,e^{i\xi_n}\bigr)$ as in (\ref{V}).   The Fourier dual operation to the pushforward $\tau_\ast$ occurring in Theorem~\ref{thm-capture} is restriction to a line with direction $\tau$; we fix notation as follows.  Given $\tau\in\real^n$, let $\ell_\tau$ denote the line on the torus defined as 
\begin{equation}\label{line}
\ell_\tau:\real\rightarrow\ttt^n,\qquad \ell_\tau(\sigma)=\bigl(e^{i\sigma\tau_1},\ldots,e^{i\sigma\tau_n}\bigr)\quad(\sigma\in\real).
\end{equation} 
The Fourier transforms $\widehat{G}_\eta$ of solutions to (\ref{model}) are obtained by restricting $\hhh$ to lines on the torus.   
\begin{thmbis}{thm-capture}
Convert an arbitrary $\eta\in\fff_n$ into a pair $(\tau,w)$ of $n$-tuples by the transformations (\ref{travel-times}) and (\ref{reflectivities}). Then $\widehat{G}_\eta(\sigma)=\hhh\bigl(w,\ell_\tau(\sigma)\bigr)$. 
\end{thmbis}
Explicitly,
\begin{equation}\label{explicit-restriction}
\widehat{G}_\eta(\sigma)=
\sum_{k\in\{1\}\times\integer^{n-1}}\Bigl(\prod_{j=1}^n\varphi^{(k_j,k_{j+1})}(w_j)\Bigr)e^{i\sigma\langle k,\tau\rangle}.
\end{equation}
This result resolves the first essential problem mentioned in \S\ref{sec-essential-questions}.  Moreover, its dual version, Theorem~\ref{thm-capture}, opens the way to a strikingly simple qualitative description of how data depends on physical parameters.

\subsection{How data depends on physical parameters\label{sec-data-dependence}}

The solution to the system (\ref{U-form}), (\ref{Helmholtz}), (\ref{coupling}), and (\ref{trace}) provides the basis for a simple geometric interpretation of the correspondence between measured data and  physical parameters.  The key idea is to encode the coefficients of $\ggg$ in Theorem~\ref{thm-unique}, defined as
\begin{equation}\label{ck}
c_k(w)=\left\{
\begin{array}{cc}
0&\mbox{ if }k_1\neq 1\\
\prod_{j=1}^n\varphi^{(k_j,k_{j+1})}(w_j)&\mbox{ if }k_1=1
\end{array}
\right.\qquad\bigl((k,w)\in\integer^n\times\overline{\ddd}^n\bigr),
\end{equation}
in a single amplitude function $\psi:\complex^n\rightarrow\complex$.  There are many ways to do this; below we illustrate one such construction based on the strict floor function.\footnote{For a different construction restricted to $\fff_n^\real$, see \cite{Gi:Dolo2014}.}   

We refer to the lattice comb
\begin{equation}\label{sampling-operator}
\mathbf{S}(z)=\sum_{k\in\integer^n}\delta(z-k)
\end{equation}
as the lattice point sampling distribution on $\complex^n$.  The term sampling distribution refers to the model for sampling whereby a pointwise well-defined function $f:\complex^n\rightarrow\complex$ is sampled on the integer lattice by way of multiplication by $\mathbf{S}$, so that the distribution $\mathbf{S}f$ represents the sampled version of $f$.  

\begin{figure}[h]
\hspace*{\fill}
\parbox{5.8in}{
\fbox{
\includegraphics[clip,trim=0in 0in 0in 0in, width=2.8in]{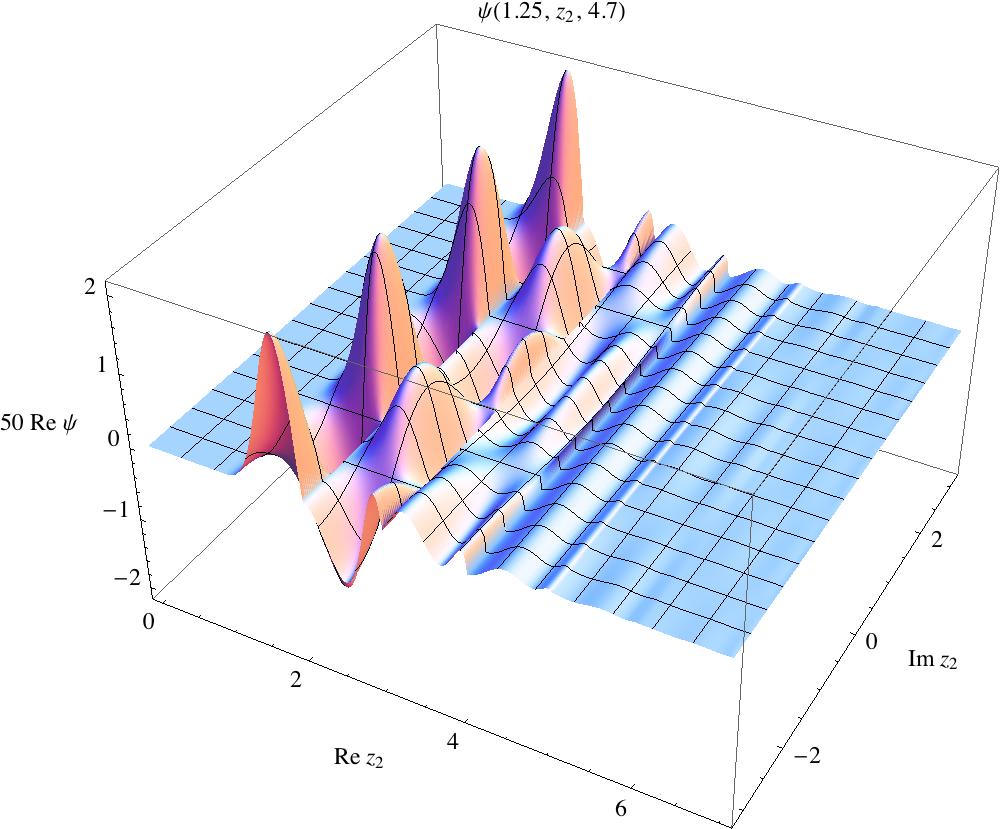}
\includegraphics[clip,trim=0in 0in 0in 0in, width=2.8in]{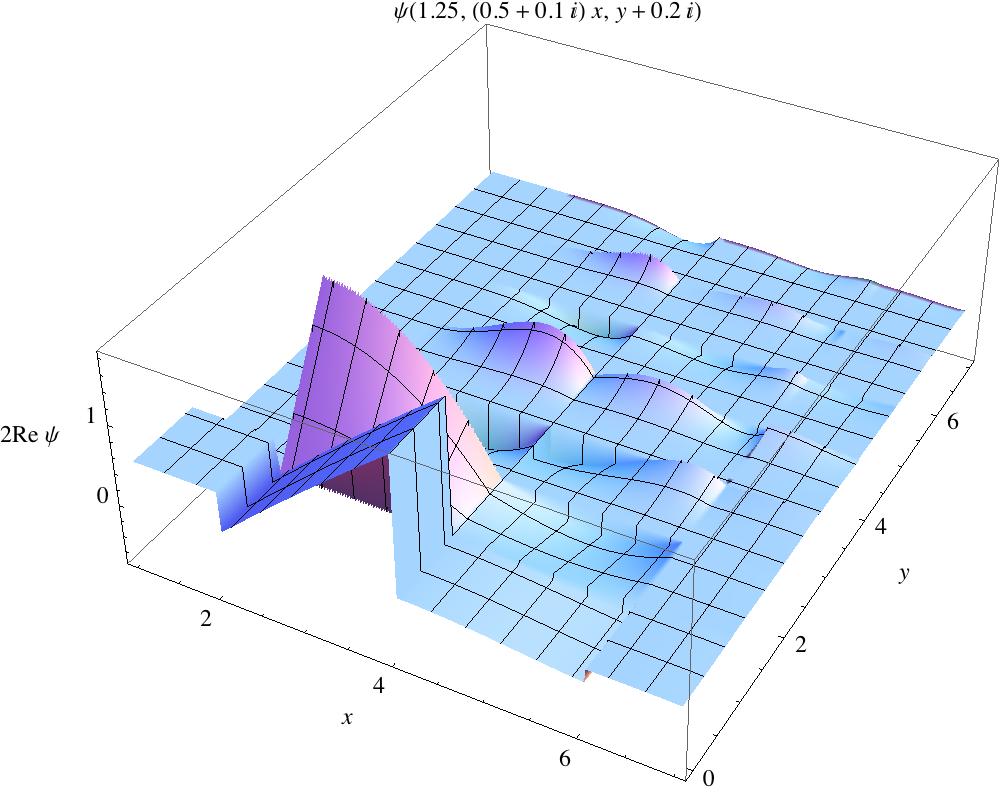}%left,bottom,right,top
}
\caption{ The amplitude wavefield on two-dimensional slices of $\complex^3$. On the left is the surface $\Re\psi(1.25,x+iy,4.7)$ for $0\leq x\leq 7$ and $-\pi\leq y\leq \pi$.  The surface plotted at right is $\Re\psi(1.25,x(.5+.1i), y+.2i)$ for $1\leq x\leq 7$ and $0\leq y\leq 7$.}\label{fig-psi}
}\hspace*{\fill}
\end{figure}

Let the symbol $\lfloor\cdot\rfloor$ denote both the strict floor function on $\real$ and its extension $\lfloor\cdot\rfloor:\real^n\rightarrow\integer^n$ defined by 
\[
\lfloor x\rfloor=\bigl(\sup\{j\in\integer\,|\,j<x_1\},\ldots,\sup\{j\in\integer\,|\,j<x_n\}\bigr).
\]
Let $w^{\circ}$ denote the polar coordinates of a point $w\in\left(\overline{\ddd}\setminus\{0\}\right)^n$ represented as a complex vector,
\begin{equation}\label{polar}
w^\circ=\bigl(|w_1|+i\arg w_1,\ldots,|w_n|+i\arg w_n\bigr),
\end{equation}
and define $\psi:\complex^n\rightarrow\overline{\ddd}$ by the formula 
\begin{equation}\label{psi_c_k}
\psi(k+w^\circ)=c_k(w)\qquad\bigl(k\in\integer^n\bigr)
\end{equation}
so that 
\begin{equation}\label{psi}
\psi(z)=c_{\lfloor\Re z\rfloor}\Bigl( \bigl(\Re z_1-\lfloor\Re z_1\rfloor\bigr)e^{i\Im z_1},\ldots,\bigl(\Re z_n-\lfloor\Re z_n\rfloor\bigr)e^{i\Im z_n}\Bigr).
\end{equation}
The explicitly defined function $\psi$ plays the role of a universal amplitude wavefield (see Figure~\ref{fig-psi}) whose translates sampled on the integer lattice push forward to $G_\eta$ for any $\eta\in\fff_n$, as follows.\footnote{The present construction of $\psi$ disallows reflectivity $w_j=0$. In terms of $\eta$ the case $w_j=0$ corresponds to $C_{j-1}=0$, which implies $\eta\in\fff_{n-1}$, i.e the case is captured by a model with fewer layers.  On the other hand, if desired, the case $w_j=0$ can be formally included by constructing a version of $\psi$ that involves a slightly more complicated change of variables---the present choice is a matter of taste.}

Letting $T_\alpha$ denote the operation of translation by $\alpha$,
\[
T_\alpha f(z)=f(z-\alpha),
\]
the solution $\ggg(w,x)$ to the system (\ref{U-form}), (\ref{Helmholtz}), (\ref{coupling}), and (\ref{trace}) may be expressed in terms of $\psi$ and the sampling distribution by the formula
\begin{equation}\label{sampled-translate}
\ggg(w,x)=\mathbf{S}T_{-w^\circ}\psi(x).
\end{equation} 
Interpreting $\tau\in\real^n$ as a linear functional on $\complex^n$ via the formula 
\[
\tau(z)=\langle \Re z,\tau\rangle,
\] 
the impulse responses $G_\eta$ is given in terms of $\psi$ by the following result.  
\begin{thm}\label{thm-psi}
Let $\psi:\complex^n\rightarrow\complex$ denote the wavefield defined explicitly by (\ref{psi}), (\ref{ck}) and (\ref{defn-scattering-polynomials}).   
Then for any $\eta\in\fff_n$, 
\[
G_\eta=\tau_\ast\textbf{S}\,T_{-w^\circ}\psi,
\]
where $\tau$ and $w$ are obtained from $\eta$ by the transformations (\ref{travel-times}) and (\ref{reflectivities}).
\end{thm}
Expanding the formula $\tau_\ast\textbf{S}\,T_{-w^\circ}\psi$ yields
\begin{equation}\label{expanded}
G_\eta(t)=\sum_{k\in\integer\rule{0pt}{5pt}^n}\psi(k+w^\circ)\delta(t-\langle k,\tau\rangle).
\end{equation}
Thus $G_\eta$ is obtained from $\psi$ by the following sequence of operations:
\begin{enumerate}
\item translation by $-w^\circ$;
\item sampling on the integer lattice;
\item pushing forward by $\tau_\ast$. 
\end{enumerate}
The three operations above comprise a qualitative description of the way the impulse response $G_\eta$ depends on the physical parameters $(\tau,w)$ that is remarkably simple.  The measured data $\chi_{[0,T]}G_\eta$ is then obtained by restricting the series (\ref{expanded}) to the finite set of indices $k$ for which 
$
0\leq\langle k,\tau\rangle\leq T.
$

\section{Proof of the main results\label{sec-proof}}

\subsection{The scattering disk\label{sec-new-polynomials}}

In the present section we give a brief overview of the scattering disk $\bigl(\overline{\ddd}, ds^2\bigr)$ with the metric (\ref{hybrid-metric}), an interesting object in its own right.  Our main goal is to prove Theorem~\ref{thm-hybrid-laplacian}; proofs of facts not directly relevant to this will be deferred to a separate paper.   

Because the scattering disk has finite diameter $2\pi$, its isometries preserve distance to the boundary, and are hence generated by rotations $z\mapsto \lambda z$ $\bigl(\lambda\in\ttt\bigr)$ and conjugation $z\mapsto\bar{z}$ as in the euclidean case.  On the other hand, the volume measure of the scattering disk,
\begin{equation}\label{area-form}
d\mu=\frac{4}{1-x^2-y^2}dxdy,
\end{equation}
is such that all even moments are infinite:
\begin{equation}\label{moments}
\int_{\ddd}x^{2m}y^{2n}\,d\mu=\infty\quad\mbox{ for every }m,n\geq 0.
\end{equation}
The theory of multivariate orthogonal polynomials is predicated on measures having finite moments, with a defining property of the set of orthogonal polynomials of total degree $n$ being orthogonality to all polynomials of degree $<n$, viz.~\cite{Ko:1975},\cite{DuXu:2001}.  Measures such as (\ref{area-form}) are excluded a priori; this serves as a plausible reason why the scattering polynomials (\ref{scattering-polynomial}) of Definition~\ref{defn-scattering-polynomials} have not been previously noticed.  
Of course, our scattering polynomials are by definition multiples of $1-x^2-y^2=1-\zeta\bar{\zeta}$, which guarantees their finite integrability with respect to $d\mu$.  

\begin{figure}[h]
\hspace*{\fill}
\parbox{5.5in}{
\fbox{
\includegraphics[clip,trim=0in 0in 0in 0in, width=2.75in]{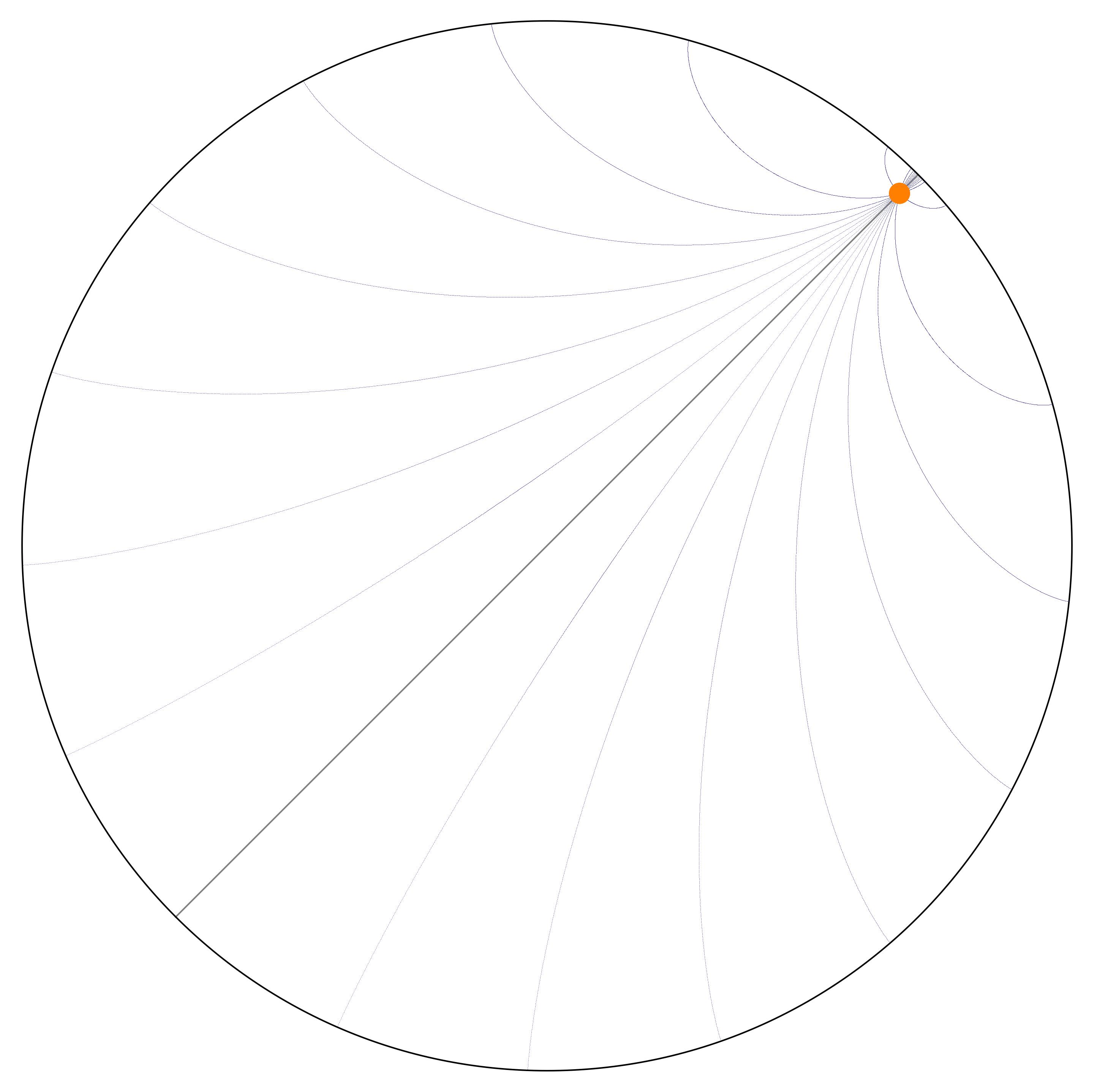}
\includegraphics[clip,trim=0in 0in 0in 0in, width=2.75in]{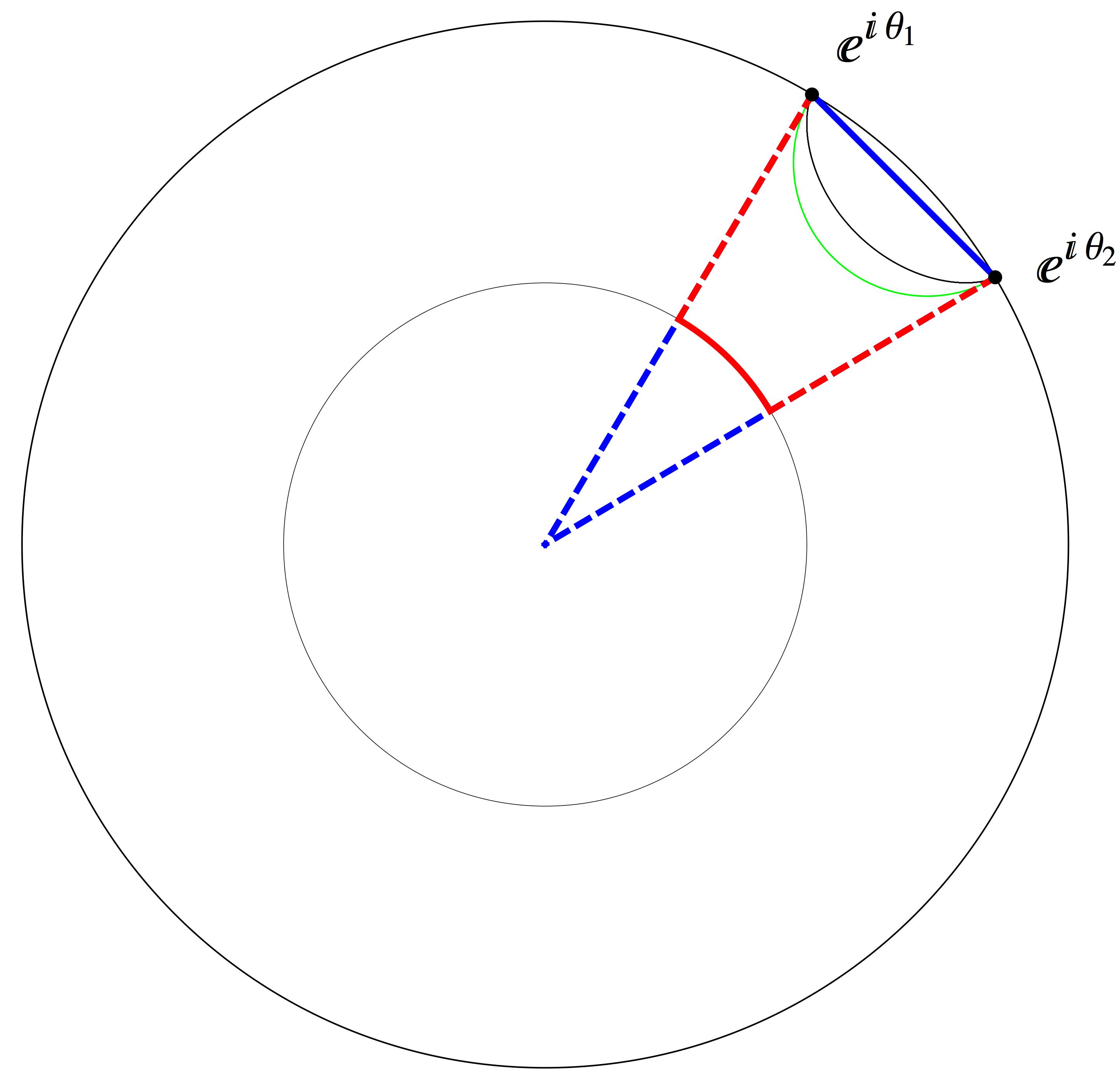}
%left,bottom,right,top
}
\caption{At left, geodesics through a point in the scattering disk.  At right, the thin black curve is the scattering geodesic connecting boundary points $e^{i\theta_1}$ and $e^{i\theta_2}$. It lies between the (green, circular) hyperbolic and the (blue) euclidean geodesics.  In the scattering metric the distance along the euclidean geodesic is the same as that along the dashed path to the centre of the disk and back; both have length $2\pi$.  Distance along the red contour gives an easy upper bound for the distance between $e^{i\theta_1}$ and 
$e^{i\theta_2}$ in the scattering metric.}\label{fig-geodesics}
}\hspace*{\fill}
\end{figure}

The scattering disk is akin to the hyperbolic disk in that both have infinite area.  The two manifolds are also similar on the level of geodesics (Fig.~\ref{fig-geodesics}).  
In terms of coordinates $\zeta=re^{i\theta}$, geodesics through a point $(r_0,\theta_0)$ in the scattering disk are given by the formula
\begin{equation}\label{geodesics}
\theta(r)=\pm\left(\arctan\left(\sqrt{\frac{c_0r^2-1}{1-r^2}}\right)-\frac{1}{\sqrt{c_0}}\arcsin\left(\sqrt{\frac{c_0r^2-1}{c_0-1}}\right)\right)+c_1\qquad\bigl(c_0\geq r_0^{-2}\bigr),
\end{equation}   
where    
\[
c_1=\theta_0\mp\left(\arctan\left(\sqrt{\frac{c_0r_0^2-1}{1-r_0^2}}\right)-\frac{1}{\sqrt{c_0}}\arcsin\left(\sqrt{\frac{c_0r_0^2-1}{c_0-1}}\right)\right).
\]
From this, one computes $\lim_{r\rightarrow1^-}\theta^\prime(r)=0$, showing that geodesics meet the boundary circle at right angles (Fig.~\ref{fig-geodesics}). 

Finally, one can compute an upper bound on the geodesic distance between nearby points on the boundary using a simple contour (see Fig.~\ref{fig-geodesics}) to conclude that scattering distance between nearby points is bounded by a constant multiple of the square root of euclidean distance.  It follows that if a neighbourhood of a boundary point is open with respect to the scattering metric, then it is open also with respect to the euclidean metric.  This in turn implies that the scattering disk is compact, unlike the hyperbolic disk.   

\subsubsection{Proof of Theorem~\ref{thm-hybrid-laplacian}\label{sec-proof-3.1}}  We turn now to the proof of Theorem~\ref{thm-hybrid-laplacian}.  Expanding the binomial $(1-\zeta\bar{\zeta})^{p+q-1}$ in formula (\ref{scattering-polynomial}) of Definition~\ref{defn-scattering-polynomials}, and then applying the derivative $\partial^{p+q}/\partial\zeta^p\partial\bar{\zeta}^q$, yields the formula
\begin{equation}\label{step-one}
\varphi^{(p,q)}(\zeta)=
\frac{(-1)^{q+\nu}}{q}(1-\zeta\bar{\zeta})\zeta^{m+\nu-p}\bar{\zeta}^{m+\nu-q}\sum_{j=0}^{\nu-1}(-1)^j\frac{(j+\nu+m)!}{j!(j+m)!(\nu-j-1)!}(\zeta\bar{\zeta})^j,
\end{equation}
where $m=|p-q|$ and $\nu=\min\{p,q\}$; the latter notation will be used in the remainder of this section.  One can verify directly that $-\mlap\varphi^{(p,q)}=pq\,\varphi^{(p,q)}$ by applying the operator 
\[
-\mlap=-(1-\zeta\bar{\zeta})\frac{\partial^2}{\partial\zeta\partial\bar{\zeta}}
\]
to the right-hand side of (\ref{step-one}).  
It also follows from (\ref{step-one}) that 
\begin{equation}\label{step-two}
\varphi^{(p,q)}\bigl(re^{i\theta}\bigr)=e^{i(q-p)\theta}f^{(p,q)}(r),
\end{equation}
where 
\begin{equation}\label{fpq}
f^{(p,q)}(r)=\frac{(-1)^{q+\nu}}{q}(1-r^2)r^m\sum_{j=0}^{\nu-1}(-1)^j\frac{(j+\nu+m)!}{j!(j+m)!(\nu-j-1)!}r^{2j}.
\end{equation}
The radial functions $f^{(p,q)}$ were first discovered in \cite{Gi:SIAP2014}, as was the following connection to Jacobi polynomials, valid for $\nu\geq 1$.   
\begin{equation}\label{Jacobi}
f^{(p,q)}(r)=\frac{(-1)^{q+\nu}(m+\nu)}{q}(1-r^2)r^mP^{(m,1)}_{\nu-1}(1-2r^2).  
\end{equation}
Note that the angular part of $\varphi^{(p,q)}(re^{i\theta})$, namely $e^{i(q-p)\theta}$, is a pure frequency.  Therefore if $q-p\neq q^\prime-p^\prime$, then $\varphi^{(p,q)}$ and $\varphi^{(p^\prime,q^\prime)}$ are orthogonal, both in $L^2\bigl(\ddd,d\mu\bigr)$ and $L^2\bigl(\ddd,dxdy\bigr)$.  In particular, if $pq=p^\prime q^\prime$ and $(p,q)\neq(p^\prime,q^\prime)$, then $\varphi^{(p,q)}$ and $\varphi^{(p^\prime,q^\prime)}$ are orthogonal, so the set of scattering polynomials corresponding to any fixed eigenvalue is linearly independent.  

In order to complete the proof of Theorem~\ref{thm-hybrid-laplacian}, it remains to show that: (1) only non-negative integers are eigenvalues of $-\mlap$ with zero boundary values; and (2) for any two integers $n\neq0$ and $k\geq 1$ there is at most one radial function $f(r)$ (up to scalar multiplication) such that $f(r)e^{in\theta}$ is an eigenvalue of $-\mlap$ with eigenvalue $k$.  Both (1) and (2) will be seen to follow from separation of variables applied to the eigenvalue equation for $-\mlap$, as follows.   

Since it is elliptic and has analytic coefficients, the operator $\mlap+k$ is analytic hypoelliptic for any constant $k$. Therefore any distributional eigenvalue $\varphi$ of $-\mlap$ is necessarily a real analytic function \cite[Thm.~10]{Ha:2006}.  Such a function is the uniform limit of its radial Fourier series
\[
\varphi(re^{i\theta})=\sum_{n\in\integer}a_n(r)e^{in\theta},
\]
on which the operator $-\mlap$ may be evaluated term by term.   The image of any particular term $a_n(r)e^{in\theta}$ by $-\mlap$ is a function $A_n(r)e^{in\theta}$ having the same angular part.  Hence if $\varphi$ is an eigenfunction of $-\mlap$ with eigenvalue $k$, the same is true of each nonzero term $a_n(r)e^{in\theta}$.  In other words, real analytic tensor products of the form $f(r)e^{in\theta}$, where $n\in\integer$, span the eigenspaces of $-\mlap$, and separation of variables is guaranteed not to miss any solutions to the eigenvalue problem.  

Suppose therefore that, for some $n\in\integer$, $f(r)e^{in\theta}$ is a real analytic eigenfunction of $-\mlap$ corresponding to non-zero eigenvalue $k$ such that $f(1)=0$.  We shall verify that the eigenvalue equation determines $f$ up to a scalar multiple, and that $k$ is necessarily a positive integer, as follows.  Expressing the hybrid laplacian in polar coordinates yields
\[
(\mlap+k)\bigl(f(r)e^{in\theta}\bigr)=\frac{1-r^2}{4}\left(f^{\prime\prime}(r)+\frac{1}{r}f^\prime(r)+\left(\frac{4k}{1-r^2}-\frac{n^2}{r^2}\right)f(r)\right)e^{in\theta}.
\]
Therefore $(\mlap+k)\bigl(f(r)e^{in\theta}\bigr)=0$ implies that $f(r)$ satisfies the equation
\begin{equation}\label{ode}
r^2(1-r^2)f^{\prime\prime}+r(1-r^2)f^\prime+\bigl(4kr^2-n^2(1-r^2)\bigr)f=0.
\end{equation}
By real analyticity $f$ has a convergent Taylor expansion,
\begin{equation}\label{taylor-expansion}
f(r)=\sum_{j=0}^\infty b_jr^j
\end{equation}
upon which equation (\ref{ode}) induces the recurrence relation
\begin{align}
n^2b_0&=0\label{b0}\\
(1-n^2)b_1&=0\label{b1}\\
\forall j\geq0\quad\bigl((j+2)^2-n^2\bigr)b_{j+2}&=\bigl(j^2-4k-n^2)b_j\label{bj+2}.
\end{align}
Assuming $f$ is not identically zero, one deduces directly from (\ref{b0}),(\ref{b1}),(\ref{bj+2}) that the least index $m$ for which $b_m\neq 0$ is $m=|n|$, with the particular choice of $b_{|n|}$ determining $f$ itself by recurrence.  This proves that there is at most one radial function $f(r)$ (up to choice of $b_{|n|}$) such that $f(r)e^{in\theta}$ is an eigenfunction of $-\mlap$ corresponding to a given non-zero eigenvalue $k$.   

Write $m=|n|$. Solving the recurrence yields $b_{m+2j+1}=0$ $(j\geq 0)$, 
\begin{equation}\label{recurrence-solution}
b_{m+2}=\frac{-k}{m+1}b_m\quad\mbox{ and }\quad 
b_{m+2j}=\frac{-k}{j(m+j)}\prod_{\nu=1}^{j-1}\left(1-\frac{k}{\nu(m+\nu)}\right)b_m\quad\quad\forall j\geq2.
\end{equation}
The boundary value $f(1)=0$ determines possible values of $k$.  By (\ref{recurrence-solution}), $f(1)=b_m\xi(k)$, where 
\begin{equation}\label{xi}
\begin{split}
\xi(k)&=1-k\left(\frac{1}{m+1}+\sum_{j=2}^\infty\frac{1}{j(m+j)}\prod_{\nu=1}^{j-1}\left(1-\frac{k}{\nu(m+\nu)}\right)\right)\\
&=\prod_{\nu=1}^\infty\left(1-\frac{k}{\nu(m+\nu)}\right).
\end{split}
\end{equation}
Convergence of the series $\sum1/(\nu(m+\nu))$ for all $\nu\geq1$ ensures that $f(1)=b_m\xi(k)=0$ if and only if $k=\nu(m+\nu)$ for some $\nu\geq1$.  This proves that $k$ is a positive integer, completing the proof of Theorem~\ref{thm-hybrid-laplacian}.  

Slightly more work produces the scattering polynomials themselves, up to a scalar multiple.  It follows from (\ref{recurrence-solution}) that 
\[
b_{m+2j}=0\mbox{ if }j>\nu\quad\mbox{ and }\quad b_{m+2j}\neq0\mbox{ if }0\leq j\leq\nu. 
\]
Thus $f(r)$ is a polynomial of precise degree $m+2\nu$ that has a zero of order $m$ at $r=0$, and such a solution exists for every pair of integers $m\geq0$ and $\nu\geq 1$.  Set
\begin{equation}\label{beta}
\beta_0=1,\qquad\beta_j=1-\frac{\nu(m+\nu)}{j(m+j)}\quad(1\leq j\leq \nu).
\end{equation}
Note in particular that $\beta_\nu=0$.   Combining the formula (\ref{recurrence-solution}) with the expansion (\ref{taylor-expansion}) shows that 
\begin{equation}\label{explicit}
f(r)=b_mr^m(1-r^2)\sum_{j=0}^{\nu-1}\Bigl(\prod_{s=0}^j\beta_s\Bigr)r^{2j},
\end{equation}
which has two associated eigenvalues $\varphi(re^{\pm im\theta})=f(r)e^{\pm im\theta}$ if $m\geq1$, and just one if $m=0$. With $m=|n|=|p-q|$ and $\nu=\min\{p,q\}$, the formula (\ref{explicit}) is proportional to $f^{(p,q)}(r)$ and $f^{(q,p)}(r)$ as defined in (\ref{fpq}).  Thus separation of variables yields the scattering polynomials (up to a scalar multiple) directly from the eigenvalue equation for $-\mlap$.

\subsection{Fourier series on the torus\label{sec-fourier-expansion}}

We shall prove Theorems~\ref{thm-unique}$^\prime$ and \ref{thm-capture}$^\prime$, since from the technical point of view it is slightly more convenient to work with the system (\ref{wave}), (\ref{coupling-dual}) and (\ref{V-boundary}) on $\overline{\ddd}^n\times\ttt^n$.  Fourier analysis on $\ttt^n$ does not involve any distinction between tempered and compactly supported distributions, since the manifold itself is compact (see \cite[Ch.~3]{Gr:2008}).  Theorems~\ref{thm-unique} and \ref{thm-capture} follow by Fourier duality.  Note that Theorem~\ref{thm-psi} follows from the latter two by construction.  

For $(w,z)\in\overline{\ddd}^n\times\ttt^n$ fix the notation
\[
w_j=|w_j|e^{i\theta_j},\quad z_j=e^{i\xi_j}\quad(1\leq j\leq n).
\]
Let $\sss_n$ denote the system (\ref{wave}), (\ref{coupling-dual}) and (\ref{V-boundary}), which we reproduce here for ease of reference:
\begin{align}
-\widetilde{\Delta}_jV+\frac{\partial^2V}{\partial\xi_{j+1}\partial\xi_j}&=0\label{proof-wave}\\
\frac{\partial V}{\partial\theta_j}+\frac{\partial V}{\partial\xi_j}-\frac{\partial V}{\partial\xi_{j+1}}&=0\label{proof-coupling}\rule{0pt}{22pt}\\
\forall w\in\partial\overline{\ddd}^n,\quad V(w,\one)&=\eee(w)\label{proof-boundary}\rule{0pt}{22pt}
\end{align}
where $1\leq j\leq n$ and $\partial/\partial\xi_{n+1}=0$.   

In the present section we make some initial observations concerning the Fourier series of the solution to $\sss_n$.  The pairing of a distribution with a test function will be denoted by round brackets, whatever the domain; we take distributions to be continuous linear functionals (as opposed to conjugate linear).   Let $V\in\dist(\overline{\ddd}^n\times\ttt^n)$ be a distribution.    The Fourier coefficients of $V$ are distributions $c_k\in\dist(\overline{\ddd}^n)$ defined by the rule
\[
(c_k,\phi)=(V,\phi\otimes\mu_k)\qquad(k\in\integer^n),
\]
where $\mu_k(z)=\overline{z}^k/(2\pi)^n=e^{-i\langle k,\xi\rangle}/(2\pi)^n$.   
In order that $V$ be considered as a candidate solution to the system $\sss_n$ it has to be possible to interpret the partial boundary condition (\ref{proof-boundary}). This requires in particular that $V$ be trace class in the sense that the series $\sum_{k\in\integer^n}c_k$
should converge weakly to a bona fide distribution in $\dist(\overline{\ddd}^n)$.  
In fact for any distributional solution $V$ to $\sss_n$ the Fourier coefficients $c_k$ must be smooth functions, as follows. In reference to lattice points $k\in\integer^n$ we use the convention established earlier whereby $k_{n+1}=0$.   

\begin{prop}\label{prop-harmonic}
Any distributional solution $V$ to the equations (\ref{proof-wave}) has real analytic Fourier coefficients $c_k$ $(k\in\integer^n)$. Each such coefficient is a tensor product univariate functions, of the form
\begin{equation}\label{coefficient-structure}
c_k(w)=\prod_{j=1}^n\phi_j(w_j),\quad\mbox{ where }\quad-\mlap_j\phi_j=k_jk_{j+1}\phi_j\quad(1\leq j\leq n),
\end{equation}
and is uniquely determined by its restriction to the set 
\[
\{w\in\overline{\ddd}^n\,|\,|w_n|=1\}.
\]
\end{prop}
\begin{pf}
If $V$ satisfies (\ref{proof-wave}) then for any $k\in\integer^n$ and any test function $\phi\in C^\infty(\overline{\ddd}^n)$
\[
\begin{split}
0&=(-\widetilde{\Delta}_jV+\frac{\partial^2V}{\partial\xi_{j+1}\partial\xi_j},\phi\otimes\mu_k)\\
&=(V,-\widetilde{\Delta}_j\phi\otimes\mu_k-k_jk_{j+1}\phi\otimes\mu_k)\\
&=(V,(-\widetilde{\Delta}_j\phi-k_jk_{j+1}\phi)\otimes\mu_k)\\
&=(c_k,-\widetilde{\Delta}_j\phi-k_jk_{j+1}\phi)\\
&=(-\widetilde{\Delta}_jc_k-k_jk_{j+1}c_k,\phi),
\end{split}
\]
which implies that 
\begin{equation}\label{coefficient-equation}
-\widetilde{\Delta}_jc_k-k_jk_{j+1}c_k=0.
\end{equation}
Since the operator $-\widetilde{\Delta}_j-k_jk_{j+1}$ is elliptic with real analytic coefficients it follows by analytic hypoellipticity that each $c_k$ is real analytic.   
Thus the boundary condition (\ref{proof-boundary}) may be interpreted in the sense of ordinary functions, where 
\[
V(w,\one)=\sum_{k\in\integer^n}c_k(w).
\]
Equations (\ref{coefficient-equation}) for $1\leq j\leq n$ imply furthermore that each $c_k$ is a tensor product of eigenfunctions of the hybrid laplacian, 
\[
c_k(w)=\prod_{j=1}^n\phi_j(w_j),
\]
with $-\widetilde{\Delta}_j\phi_j=k_jk_{j+1}\phi_j$.   Note in particular that $k_{n+1}=0$.  Therefore $\phi_n$ is a harmonic function, determined by its restriction to the boundary circle $|w_n|=1$.   
\end{pf}

The real analytic functions $c_k$ may themselves be expanded as Fourier series with radial coefficients in the form
\begin{equation}\label{coefficient-expansion}
c_k(w)=\sum_{l\in\integer^n}d_{k,l}(r)e^{i\langle l,\theta\rangle},
\end{equation}
where $r=(|w_1|,\ldots,|w_n|)$.  Denote the left shift operator by a tilde, so that 
\[
\tilde{k}=(k_2,k_3,\ldots,k_n,0).
\]
\begin{prop}\label{prop-V-coupling}
Equations (\ref{proof-coupling}) imply that in the expansion (\ref{coefficient-expansion}) $d_{k,l}=0$ unless $l=\tilde{k}-k$.  
\end{prop}
\begin{pf}
For $w\in\overline{\ddd}^n$ write 
\[
r=(|w_1|,\ldots,|w_n|),\qquad v=\left(\frac{w_1}{|w_1|},\ldots,\frac{w_n}{|w_n|}\right),
\]
and let $P:\overline{\ddd}^n\rightarrow[0,1]^n\times\ttt^n$ denote the change of variables
\[
w\stackrel{P}{\mapsto}(r,v).
\]
The coefficient $d_{k,l}$ is defined in terms of $V$ by the formula
\begin{equation}\label{d-formula}
(d_{k,l},\rho)=\bigl(c_k,(\rho\otimes\mu_l)\circ P\bigr)=\bigl(V,((\rho\otimes\mu_l)\circ P)\otimes\mu_k\bigr)
\end{equation}
for test functions $\rho\in C^\infty([0,1]^n)$.  If $V$ satisfies the equation (\ref{proof-coupling}) it follows that 
\[
\begin{split}
0&=(C_jV,((\rho\otimes\mu_l)\circ P)\otimes\mu_k)\\
&=(V,-C_j((\rho\otimes\mu_l)\circ P)\otimes\mu_k)\\
&=(V,-(l_j+k_j-k_{j+1})((\rho\otimes\mu_l)\circ P)\otimes\mu_k)\\
&=(d_{k,l},-(l_j+k_j-k_{j+1})\rho)\\
&=(-(l_j+k_j-k_{j+1})d_{k,l},\rho),
\end{split}
\]
where 
\[
C_j=\frac{\partial }{\partial\theta_j}+\frac{\partial }{\partial\xi_j}-\frac{\partial }{\partial\xi_{j+1}}.
\]
Therefore $-(l_j+k_j-k_{j+1})d_{k,l}=0$ for $1\leq j\leq n$, so that $d_{k,l}=0$ unless $l=\tilde{k}-k$, as claimed.  \end{pf}

Relabeling appropriately, we may thus express the functions $c_k$ in the form 
\begin{equation}\label{d_k}
c_k(w)=d_k(r)e^{i\langle\tilde{k}-k,\theta\rangle}=\prod_{j=1}^n\phi_j(w_j),
\end{equation}
consistent with (\ref{coefficient-structure}).  This yields a formal Fourier series for $V$ having smooth coefficients,
\begin{equation}\label{formal-fourier}
V(w,z)\sim\sum_{k\in\integer}d_k(r)e^{i\langle\tilde{k}-k,\theta\rangle}z^k.
\end{equation}
The partial boundary condition (\ref{proof-boundary}) requires that for $w\in\partial\overline{\ddd}^n$, 
\begin{equation}\label{boundary-revised}
\sum_{k\in\integer}d_k(r)e^{i\langle\tilde{k}-k,\theta\rangle}=\eee(w).
\end{equation}
In conjunction with Theorem~\ref{thm-hybrid-laplacian}, equation (\ref{d_k}) shows that if $k_jk_{j+1}\geq 1$, then $\phi_j$ has angular part $e^{i(k_{j+1}-k_j)\theta_j}$ and is hence proportional to $\varphi^{(k_j,k_{j+1})}$.  Also, if $k_j>0$ and $k_{j+1}=0$, then $\phi_j$ is proportional to $\varphi^{(k_j,0)}$, since, up to a scalar multiple, the latter is the unique harmonic function with angular part $e^{-ik_j\theta_j}$.  Lastly, if $k_jk_{j+1}<0$ for any $1\leq j\leq n$, then $\phi_j=0$ and $c_k=0$.   Thus $c_k\neq0$ only for lattice points $k$ for which $k_jk_{j+1}$ has constant sign.   

The structure of the layer collapse function $\eee$ further restricts the set of lattice points $k$ at which $c_k\neq0$ through condition (\ref{boundary-revised}).  This is worked out in detail in the next section, following a proof that the system $\sss_n$ has a unique distributional solution.

\subsection{The layer collapse function\label{sec-layer-collapse-function}}

Recall from (\ref{Psi}) and (\ref{eee-definition}) that for $w\in\overline{\ddd}^n$,
\[
\eee(w)=\Psi^{w_1}_1\circ\Psi^{w_2}_1\circ\cdots\circ\Psi^{w_n}_1(0),\quad\mbox{ where }\quad
\Psi_{1}^{w_j}(v)=\frac{v+\overline{w}_j}{1+w_jv}\quad(1\leq j\leq n).
\]
Since $\Psi^{w_n}_1(0)=\overline{w}_n$, the layer collapse function $\eee(w)$ is antiholomorphic in  the variable $w_n$, and $\eee(w)$ is therefore determined by its restriction to the set 
\[
\{w\in\overline{\ddd}^n\,|\,|w_n|=1\}.
\]
Proposition~\ref{prop-harmonic} then implies that conditions (\ref{proof-boundary}) and (\ref{boundary-revised}) extend to all of $\overline{\ddd}^n$, leading to the following uniqueness result for the system $\sss_n$.  
\begin{prop}\label{prop-uniqueness}
The system $\sss_n$ has at most one distributional solution, necessarily a function 
\[
V\in L^2(\overline{\ddd}^n\times\ttt^n)
\]
of the form
\[
V(w,z)=\sum_{k\in\integer^n}d_k(r)e^{i\langle\tilde{k}-k,\theta\rangle}z^k,
\]
where 
\[
d_k(r)=\frac{1}{(2\pi)^n}\int_{\ttt^n}\eee\bigl(r_1e^{i\theta_1},\ldots,r_ne^{i\theta_n}\bigr)e^{-i\langle\tilde{k}-k,\theta\rangle}\,d\theta.
\]
\end{prop}
\begin{pf}
Note that for every $w\in\overline{\ddd}^n$, $\eee(w)$ is a composition of disk automorphisms, evaluated at $0$, so that $|\eee(w)|\leq 1$.   
Since the coefficients $d_{k,l}$ are uniquely determined by the formula (\ref{d-formula}), the condition (\ref{boundary-revised})---extended to $\overline{\ddd}^n$---yields the given integral formula.   That $V\in L^2(\overline{\ddd}^n\times\ttt^n)$ then follows from the bound $|\eee(w)|\leq 1$.  In detail,
$(2\pi)^n\geq||\eee(re^{i\cdot})||^2_{L^2(\ttt^n)}=(2\pi)^n\sum_{k\in\integer}|d_k(r)|^2$.  Therefore,
\[
\begin{split}
||V(w,z)||^2_{L^2(\ddd^n\times\ttt^n)}&=\int_{\ddd^n\times\ttt^n}\bigl|\sum c_k(w)z^k\bigr|^2\,dwdz=(2\pi)^n\int_{\ddd^n}\sum|c_k(w)|^2\,dw\\
&=(2\pi)^{2n}\int_{[0,1]^n}\sum|d_k(r)|^2\,r_1\cdots r_n\,dr\leq 2^n\pi^{2n}.
\end{split}
\]
\end{pf}

The next step is to analyze the structure of $\eee$ to determine the set of lattice points $k\in\integer^n$ for which $c_k\neq0$.   
\vspace*{5pt}
Let $\lll_n$ denote the set of all lattice points $k\in\integer^n$ with the properties: 
\begin{enumerate}
\item each $k_j\geq0$;
\item  $k_1=1$;
\item for each $2\leq j\leq n-1$, if $k_j=0$ them $k_{j+1}=0$.
\end{enumerate}  
For $w=\bigl(r_1e^{i\theta_1},\ldots,r_ne^{i\theta_n}\bigr)\in\overline{\ddd}^n$, let $g_k(r)$ denote the $k$th Fourier coefficient of 
\[
\eee(w)=\sum_{k\in\integer^n}g_k(r)e^{i\langle k,\theta\rangle}.
\]
\begin{prop}\label{prop-index-form}
If $g_l\neq0$ then there is a unique lattice point $k\in\lll_n$ such that $l=\widetilde{k}-k$.   
\end{prop}
\begin{pf}
If $l=\tilde{k}-k$ then the conditions that characterize $k\in\lll_n$ translate in terms of $l$ to:
\begin{enumerate}
\item for each $1\leq j\leq n$, $l_1+\cdots+l_j\geq -1$;
\item $l_1+\cdots+l_n=-1$;
\item for each $1\leq j\leq n-1$, if $l_1+\cdots+l_j=-1$ then $l_1+\cdots+l_{j+1}=-1$.  
\end{enumerate}
We argue by induction on the sequence of functions
\[
v=\Psi^{w_j}_{1}\circ\cdots\circ\Psi^{w_n}_{1}(0)
\]
as $j$ decreases from $n$ to $1$.  In the base case $v=\Psi^{w_n}_{1}(0)=\overline{w_n}=e^{i(-1)\theta_n}$ there is a single one-dimensional vector, $(-1)$, corresponding to a non-zero coefficient, and the associated set $\{(-1)\}$ conforms the prescribed criteria.  For the induction, the key observation on the level of formulas is simply that 
\[
\begin{split}
\frac{\overline{w_j}+v}{1+w_jv}&=\overline{w_j}+(1-r_j^2)\frac{v}{1+w_jv}\\
&=\overline{w_j}+(1-r_j^2)(v-w_jv^2+w_j^2v^3-w_j^3v^4+w_j^4v^5-\cdots).
\end{split}
\]
Each term $w_j^sv^{s+1}$ corresponds to indices $(l_j,\ldots,l_n)$ where $l_j=s$.  By inductive hypothesis we assume that conditions (1.-3.) are satisfied for the indices $(l^\prime_{j+1},\ldots,l^\prime_n)$ occurring in $v$.  It follows by the above formula that they are again satisfied for $(l_j,\ldots,l_n)$.  For example, the total sum being -1, means that the part of the sum coming from $v$ in a term $w_j^sv^{s+1}$ is $-s-1$, to which is added $+s$ from the term $l_j=s$, for a total of $-1$.  The other items are similar.   
\end{pf}

By equation (\ref{d_k}), Proposition~\ref{prop-index-form} implies that $c_k\neq0$ only if $k\in\lll_n$.   For such lattice points every function $\phi_j(w_j)$ in (\ref{d_k}) is proportional to a scattering polynomial, since only nonnegative indices $k_j$ occur, and the case $k_j=0, k_{j+1}>0$ is ruled out.  We thus have,
\begin{cor}\label{cor-index-form}
The system $\sss_n$ has at most one distributional solution, necessarily of the form
\begin{equation}\label{almost-structure}
V(w,z)=\sum_{k\in\lll_n}\beta_k\Bigl(\prod_{j=1}^n\varphi^{(k_j,k_{j+1})}(w_j)\Bigr)z^k,
\end{equation}
where each $\beta_k$ is a scalar determined by the equation 
\begin{equation}\label{beta_k-computation}
\beta_k\prod_{j=1}^n\varphi^{(k_j,k_{j+1})}(w_j)=\frac{1}{(2\pi)^n}\int_{\ttt^n}\eee\bigl(r_1e^{i\theta_1},\ldots,r_ne^{i\theta_n}\bigr)e^{-i\langle\tilde{k}-k,\theta\rangle}\,d\theta.
\end{equation}
\end{cor}
Indeed the weights $\alpha_{p,q}$ in Definition~\ref{defn-scattering-polynomials} are chosen so that $\beta_k=1$ for every $k\in\lll_n$.  Also, since $\varphi^{(k_j,k_{j+1})}=0$ by definition if either $\min\{k_j,k_{j+1}\}<0$ or both $k_j=0$ and $k_{j+1}>0$, Corollary~\ref{cor-index-form} remains true if in (\ref{almost-structure}) the index set $\lll_n$ is replaced by the set $\{1\}\times\integer^{n-1}$, as appears in the statements of Theorems~\ref{thm-unique} and \ref{thm-capture}.    

Thus far it has been proven that $\sss_n$ has a unique distributional solution, that this solution is an $L^2$ function with smooth Fourier coefficients, and that the Fourier coefficients are tensor products of scattering polynomials.  It remains to show how this solution relates to the original model (\ref{model}).  In the next section we define a function $\kkk(w,z)$ on $\overline{\ddd}^n\times\ttt^n$ essentially by variation of parameters in the complexified backward recurrence formula.  We verify directly that $\kkk$ solves $\sss_n$, and use uniqueness to conclude that $\kkk$ coincides with (\ref{almost-structure}).   It is then a trivial matter to see that $\kkk$ restricts to the complexified backward recurrence and recovers the Fourier transforms of the solutions to the original model.  In the process, we get an easy proof that $\beta_k=1$, as claimed above, by appealing to known results.  (In principle, $\beta_k=1$ may be determined directly from (\ref{beta_k-computation}).)

\subsection{A generalized backward recurrence formula\label{sec-generalized-recurrence}}

Recall the formula (\ref{Psi}), whereby for each pair $(w,z)\in\overline{\ddd}^n\times\ttt^n$, and each $v\in\overline{\ddd}$, we write 
\[
\Psi_{z_j}^{w_j}(v)=z_j\frac{v+\overline{w}_j}{1+w_jv}\quad\quad(1\leq j\leq n).
\]
In terms of this notation define $\kkk:\overline{\ddd}^n\times\ttt^n\rightarrow\overline{\ddd}$ by
\begin{equation}\label{K}
\kkk(w,z)=\Psi_{z_1}^{w_1}\circ\Psi_{z_2}^{w_2}\circ\cdots\circ\Psi_{z_n}^{w_n}(0).
\end{equation}
Thus $\kkk$ is a variant of the backward recurrence formula (\ref{backward}) in which the terms $e^{i\sigma\tau_1},\ldots,e^{i\sigma\tau_n}$, corresponding to the line on the torus $\ell_\tau(\sigma)$, are replaced by arbitrary parameters $z_1,\ldots,z_n\in\ttt^n$.  
We shall verify that $\kkk$ is a solution to $\sss_n$, fixing notation for the relevant operators as
\[
\begin{split}
L_j&=-\widetilde{\Delta}_j+\frac{\partial^2}{\partial\xi_{j+1}\partial\xi_j}\\
C_j&=\frac{\partial }{\partial\theta_j}+\frac{\partial }{\partial\xi_j}-\frac{\partial }{\partial\xi_{j+1}}
\end{split}
\]
where $1\leq j\leq n$.

To begin we introduce a nonlinear first-order operator that will play an important auxiliary role in our analysis.  For each $1\leq j\leq n$, let $\nabla_j$ denote the gradient operator with respect to $w_j=x_j+iy_j\in\ddd$, so that 
\[
\nabla_jv=\left(\frac{\partial v}{\partial x_j},\frac{\partial v}{\partial y_j}\right)
\]
and define the differential operator $E_j$ on distributions on $\overline{\ddd}^n\times\ttt$ by the formula
\[
E_jv=-\frac{1-r_j^2}{4}\left(\nabla_jv\cdot\nabla_jv\right)+\frac{\partial v}{\partial\xi_{j+1}}\frac{\partial v}{\partial \xi_j}.
\]
(This is essentially the quadratic form associated with $L_j$, in the sense that $\langle Ljv,v\rangle=\int E_jv$.)   
\begin{lem}\label{lem-E_j}
$E_j\kkk=0$ for every $1\leq j\leq n$.  
\end{lem}
\begin{pf}
We argue by induction on the sequence of functions
\[
v=\Psi^{w_j}_{z_j}\circ\cdots\circ\Psi^{w_n}_{z_n}(0)
\]
as $j$ decreases from $n$ to $1$.  The base case $v=\Psi^{w_n}_{z_n}(0)=z_n\overline{w_n}$ trivially satisfies $E_jv=0$ if $j<n$.   Since $\frac{\partial v}{\partial\xi{n+1}}=0$ (the variable $\xi_{n+1}$ being nonexistent), for $j=n$ it suffices to verify that 
\[
-\frac{1-r_n^2}{4}\left(\nabla_nv\cdot\nabla_nv\right)=0,
\]
which follows from the fact that $\nabla_nv=(z_n,-iz_n)$.   

Next suppose that $1<j+1\leq n$ and that $E_sv^\prime=0$ for every $1\leq s\leq n$, where 
\[
v^\prime=\Psi^{w_{j+1}}_{z_{j+1}}\circ\cdots\circ\Psi^{w_n}_{z_n}(0).
\]
We will show that necessarily $E_sv=0$ where 
\[
v=\Psi^{w_j}_{z_j}\circ\cdots\circ\Psi^{w_n}_{z_n}(0).
\]
We may assume that $s\geq j$, since otherwise there is nothing to prove. 

 Suppose first that $s>j$.  By definition, 
\[
v=\Psi_{z_j}^{w_j}(v^\prime)=z_j\frac{\overline{w_j}+v^\prime}{1+w_jv^\prime}.
\]
It follows by direct computation that
\[
-\frac{1-r_s^2}{4}\left(\nabla_sv\cdot\nabla_sv\right)=
z_j^2\frac{(1-|w_j|^2)^2}{(1+w_jv^\prime)^4}\frac{-(1-r_s^2)}{4}\left(\nabla_sv^\prime\cdot\nabla_sv^\prime\right),
\]
while
\[
\frac{\partial v}{\partial\xi_{s+1}}\frac{\partial v}{\partial\xi_{s}}\frac{\partial v}{\partial \xi_s}=z_j^2\frac{(1-|w_j|^2)^2}{(1+w_jv^\prime)^4}\frac{\partial v^\prime}{\partial\xi_{s+1}}\frac{\partial v^\prime}{\partial\xi_{s}}.
\]
The induction hypothesis therefore guarantees that $E_sv=0$.   

It remains to consider the case $s=j$.  Noting that $r_j=|w_j|$, straightforward computation yields that 
\[
-\frac{1-r_j^2}{4}\left(\nabla_jv\cdot\nabla_jv\right)=
\frac{(1-r_j^2)z_j^2(\overline{w_j}+v^\prime)v^\prime}{(1+w_jv^\prime)^3},
\]
while
\[
\frac{\partial v}{\partial\xi_{j+1}}\frac{\partial v}{\partial\xi_{j}}=\frac{-z_j^2(1-r_j^2)(\overline{w_j}+v^\prime)v^\prime}{(1+w_jv^\prime)^3},
\]
from which follows the desired result that $E_jv=0$. 
\end{pf}

\begin{prop}\label{prop-L_j}
$L_j\kkk=0$ for every $1\leq j\leq n$.
\end{prop}
\begin{pf}
This is proved by an induction along the same lines as Lemma~\ref{lem-E_j}, on the sequence of functions
\[
v=\Psi^{w_j}_{z_j}\circ\cdots\circ\Psi^{w_n}_{z_n}(0)
\]
as $j$ decreases from $n$ to $1$.  The base case $v=\Psi^{w_n}_{z_n}(0)=z_n\overline{w_n}$ trivially satisfies $L_jv=0$ if $j<n$.   Since $\frac{\partial v}{\partial\xi{n+1}}=0$, for $j=n$ it suffices to verify that 
\[
-\widetilde{\Delta}_j=0,
\]
which follows from the fact that $z_n\overline{w_n}$ is harmonic with respect to $w_n$.   

Next suppose that $1<j+1\leq n$ and that $L_sv^\prime=0$ for every $1\leq s\leq n$, where 
\[
v^\prime=\Psi^{w_{j+1}}_{z_{j+1}}\circ\cdots\circ\Psi^{w_n}_{z_n}(0).
\]
We will show that necessarily $L_sv=0$ where 
\[
v=\Psi^{w_j}_{z_j}\circ\cdots\circ\Psi^{w_n}_{z_n}(0).
\]
We may assume that $s\geq j$, since otherwise there is nothing to prove. 

 Suppose first that $s>j$.  Using that  
\[
v=\Psi_{z_j}^{w_j}(v^\prime)=z_j\frac{\overline{w_j}+v^\prime}{1+w_jv^\prime},
\]
direct computation shows
\[
-\widetilde{\Delta}_sv=z_j(1-|w_j|^2)\left(\frac{-\widetilde{\Delta}_sv^\prime}{(1+w_jv^\prime)^2}-\frac{-2w_j\frac{(1-r_s^2)}{4}\left(\nabla_sv^\prime\cdot\nabla_sv^\prime\right)}{(1+w_jv^\prime)^3}\right),
\]
while
\[
\frac{\partial^2v}{\partial\xi_{s+1}\partial\xi_s}=z_j(1-|w_j|^2)\left(\frac{\frac{\partial^2v}{\partial\xi_{s+1}\partial\xi_s}}{(1+w_jv^\prime)^2}-\frac{2w_j\frac{\partial v^\prime}{\partial\xi_{s+1}}\frac{\partial v^\prime}{\partial\xi_{s}}}{(1+w_jv^\prime)^3}\right).
\]
Adding the above two parts, the induction hypothesis then implies
\[
L_sv=\frac{-2z_jw_j(1-|w_j|^2)}{(1+w_jv^\prime)^3}E_sv^\prime,
\]
which is 0 by Lemma~\ref{lem-E_j}.   

It remains to consider the case $s=j$.  Straightforward computation yields
\[
-\widetilde{\Delta}_jv=\frac{(1-|w_j|^2)z_jv^\prime}{(1+w_jv^\prime)},
\]
while
\[
\frac{\partial^2v}{\partial\xi_{j+1}\partial\xi_j}=-\frac{(1-|w_j|^2)z_jv^\prime}{(1+w_jv^\prime)^2},
\]
whence $L_jv=0$, completing the proof. \end{pf}

\begin{prop}\label{prop-C_j}
$\ccc_j\kkk=0$ for every $1\leq j\leq n$.
\end{prop}
\begin{pf}
As in the foregoing results, this may be proved by downward induction on 
\[
v=\Psi^{w_j}_{z_j}\circ\cdots\circ\Psi^{w_n}_{z_n}(0).
\]
The appropriate induction hypothesis in this case is that $C_sv=0$ for every $s\geq j$.  It then follows by direct computation that $C_s\Psi_{z_{j-1}}^{w_{j-1}}(v)=0$ for every $s\geq j-1$, yielding the desired result.  (The argument is trivial except in the case $s=j-1$, which is straightforward to check.)
\end{pf}

Lastly, note that $\kkk(w,\one)=\eee(w)$.  Thus $\kkk$ satisfies the partial boundary condition (\ref{proof-boundary}) and hence, by Propositions~\ref{prop-L_j} and \ref{prop-C_j}, the full system $\sss_n$.   It follows by Corollary~\ref{cor-index-form} that 
\begin{equation}\label{K-expansion}
\kkk(w,z)=\sum_{k\in\{1\}\times\integer^{n-1}}\beta_k\Bigl(\prod_{j=1}^n\varphi^{(k_j,k_{j+1})}(w_j)\Bigr)z^k.
\end{equation}
In particular, for real-valued $r\in[0,1]^n$ and $\tau\in\real_+^n$, 
\begin{equation}\label{real-K-expansion}
\kkk\bigl(r,(e^{i\sigma\tau_1},\ldots,e^{i\sigma\tau_n})\bigr)=\sum_{k\in\{1\}\times\integer^{n-1}}\beta_k\Bigl(\prod_{j=1}^nf^{(k_j,k_{j+1})}(r_j)\Bigr)e^{i\langle k,\tau\rangle\sigma},
\end{equation}
where the radial functions $f^{(k_j,k_{j+1})}$ are defined in \S\ref{sec-proof-3.1} by (\ref{fpq}).  
The left-hand side of (\ref{real-K-expansion}) is a particular case of the backward recurrence formula (\ref{backward}) for $\widehat{G}_\eta(\sigma)$, where $\eta$ and $(\tau,r)$ are related by the transformations (\ref{travel-times}) and (\ref{reflectivities}).  It is known \cite[Thm.~4.3]{Gi:SIAP2014} that for real-valued $r\in[0,1]^n$ and $\tau\in\real_+^n$
\begin{equation}\label{real-solution}
\widehat{G}_\eta(\sigma)=\sum_{k\in\{1\}\times\integer^{n-1}}\Bigl(\prod_{j=1}^nf^{(k_j,k_{j+1})}(r_j)\Bigr)e^{i\langle k,\tau\rangle\sigma}.
\end{equation}
Linear independence of the exponentials $e^{i\langle k,\tau\rangle\sigma}$ $(\sigma\in\real)$ for appropriately chosen $\tau$ and uniqueness of the solution to $\sss_n$ then combine to force $\beta_k=1$ for every $k\in\{1\}\times\integer^{n-1}$.  Thus 
\begin{equation}\label{finally}
\kkk(w,z)=\hhh(w,z)=\sum_{k\in\{1\}\times\integer^{n-1}}\Bigl(\prod_{j=1}^n\varphi^{(k_j,k_{j+1})}(w_j)\Bigr)z^k,
\end{equation}
completing the proof of Theorem~\ref{thm-unique}$^\prime$.   

Theorem~\ref{thm-capture}$^\prime$ then follows from the observation that for any $w\in\overline{\ddd}$ and any $\tau\in\real_+^n$, 
\[
\widehat{G}_\eta(\sigma)=\kkk\bigl(w,(e^{i\sigma\tau_1},\ldots,e^{i\sigma\tau_n})\bigr),
\]
where $\eta$ and $(\tau,w)$ are related by the transformations (\ref{travel-times}) and (\ref{reflectivities}).

\section{Conclusions\label{sec-conclusions}}

The results in the present paper have several aspects.  Firstly, they provide tools, in the guise of explicit formulas, for application to the inverse problem for the classical model (\ref{model}) for waves in layered media.  Secondly, they open a new perspective on the question of how data depends on coefficients, showing that the dependency can be described in terms of the conceptually simple---and mathematically familiar---operations of translation and pushforward (or dually, restriction), and that the nonlinear dependence is itself governed by a system of linear equations.  Thirdly, the results are illuminating from the perspective of PDEs having discontinuous coefficients.  The existence of a smooth global model is in effect a resolution of singularities; the global model is non-singular and solvable by ordinary separation of variables, and its solutions push forward (or, in the dual domain, restrict) to those of the singular model.  Lastly, the previously unknown scattering disk turns out to have an essential role underlying scattering in layered media.  In the present section we elaborate briefly on these various aspects of the paper's main results.  

\subsection{Explicit formulas and the origins of irregular structure\label{sec-origins}}

The explicit formulas for $G_\eta(t)$ and $\widehat{G}_\eta(\sigma)$ following from Theorems~\ref{thm-capture} and \ref{thm-capture}$^\prime$ are useful as a means to solve the inverse problem of determining $\eta$ from data of the form $\chi_{[0,T]}G_\eta$.  Indeed, by exploiting the connection (\ref{Jacobi}) to Jacobi polynomials, one can 
compute $\eta$ in terms of delayed data of the form $\chi_{[T_1,T_2]}G_\eta$ with $T_1>0$, which was not previously known to be possible.  (Details of this will be treated in a separate paper.)  

A basic feature of the formula (\ref{explicit-pushforward}), and its dual (\ref{explicit-restriction}), is that the singular support of $\eta$, which transforms by (\ref{travel-times}) into travel times $\tau$, and the values of $\eta$, which transform by (\ref{reflectivities}) to reflectivities, affect the data in two independent ways.  The travel time vector $\tau$ determines the support of $G_\eta$, or equivalently, the almost periods of $\widehat{G}_\eta$, while the reflectivities $w$ determine associated amplitudes.

\begin{figure}[h]
\hspace*{\fill}
\parbox{6in}{
\fbox{
\includegraphics[clip,trim=7in 6in 2.5in 3in, width=1.6in]{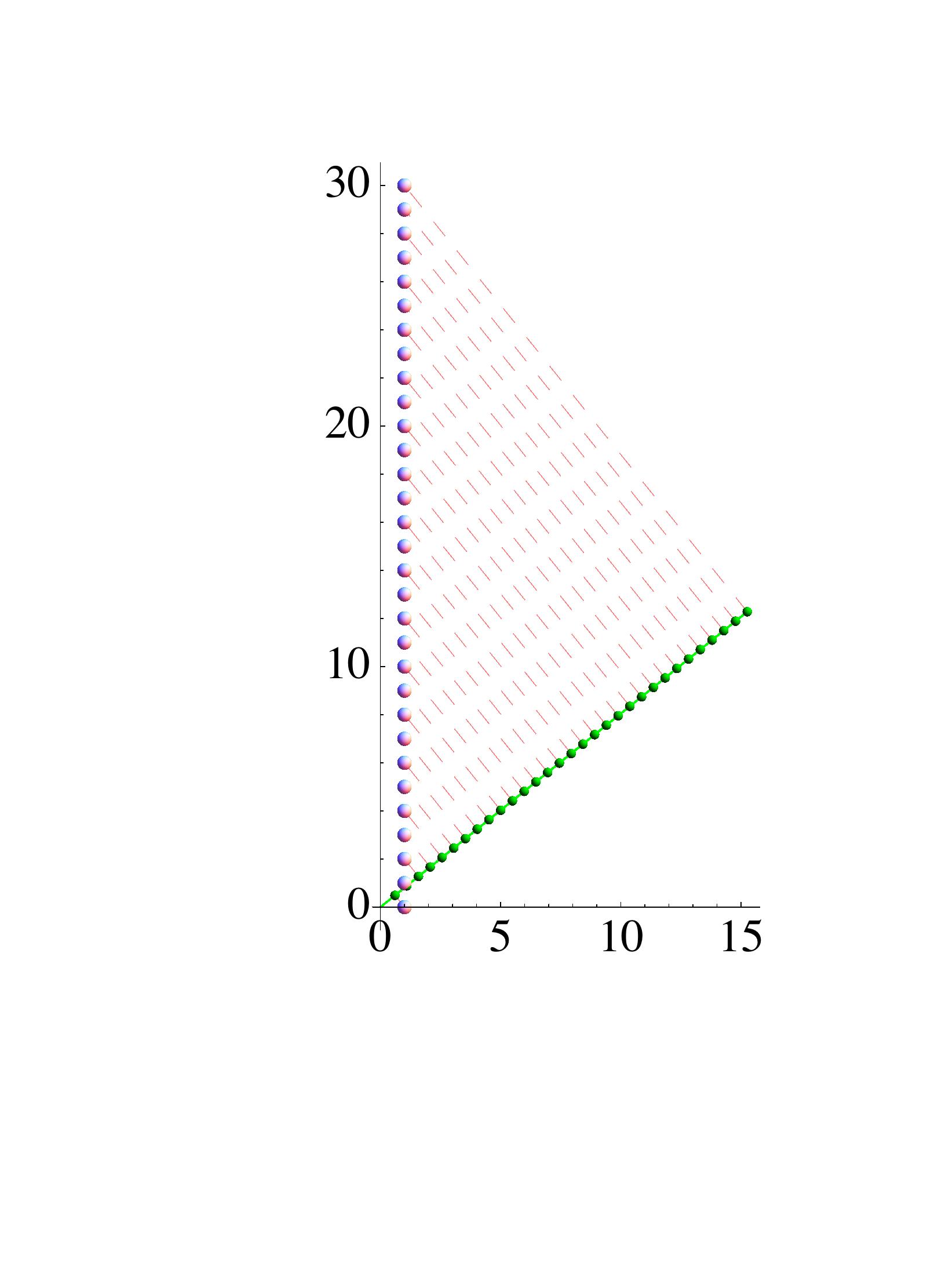}
\includegraphics[clip,trim=5in 5in 0in 18in, width=1.6in]{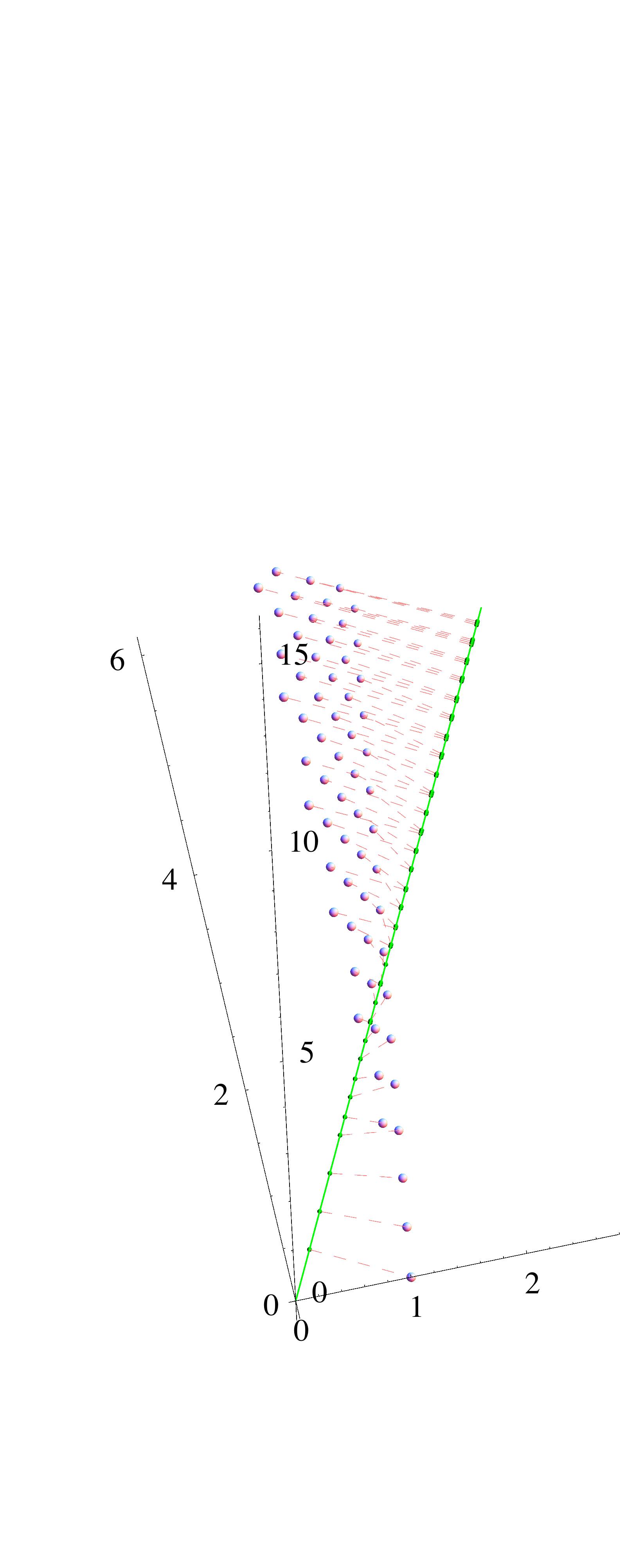}
\includegraphics[clip,trim=0in 3in 18in 28in, width=2.8in]{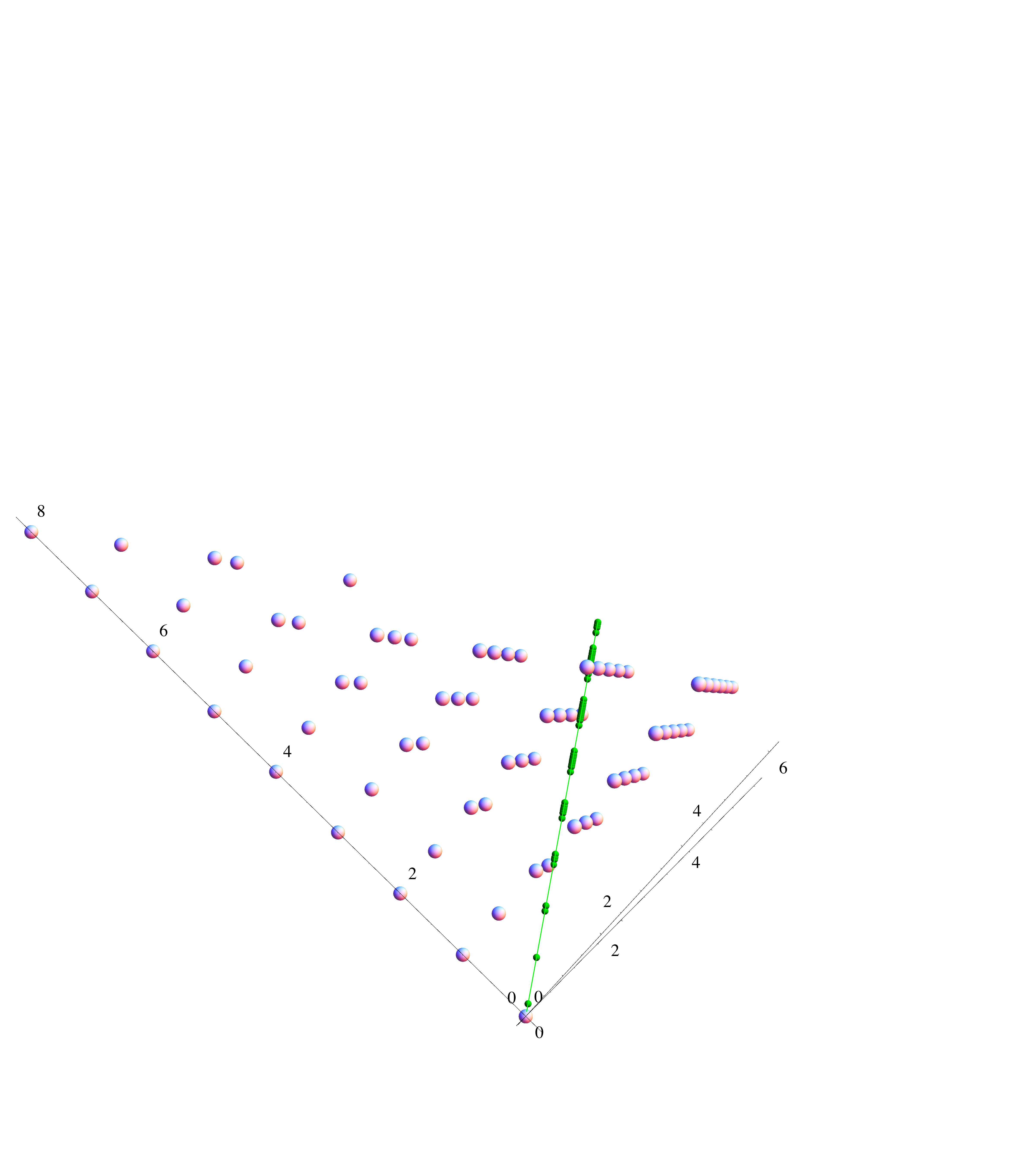}
%left,bottom,right,top
}
\caption{Examples of the pushforward of lattice points in $\lll_n$ (coloured pink) onto a timeline. The (green) image points are plotted on a line with direction $\tau\in\real^n$ for $n=2,3,4$ (from left to right), which shows the correct relative spacing. In dimension $n=2$ the spacing is regular irrespective of $\tau$; in higher dimensions it becomes irregular.}\label{fig-projections}
}\hspace*{\fill}
\end{figure}

The irregular spacing of the support of $G_\eta(t)$ is a one-dimensional manifestation---via the pushforward $\tau_\ast$---of a regular structure, namely the set of lattice points $\lll_n$, existing in higher dimensions; see Figure~\ref{fig-projections}.  Dually, the irregular almost periodicity of $\widehat{G}_\eta(\sigma)$ stems from the restriction of a regular higher dimensional structure, namely the $\xi$-periodic function $\widehat{G}(w,\xi)$, to a line $\xi=\sigma\tau$.    

\subsection{Non-linear dependence of data on coefficients\label{sec-nonlinear}}

Theorem~\ref{thm-psi} describes the data as a function of coefficients in terms of the precise, conceptually simple notions of translation and pushforward, whereas existing literature on inverse problems typically contents itself with a nondescript general designation of nonlinearity.  

Our results indicate that the values of $\eta$ affect the function $G_\eta(t)$ in a way similar to the independent variable $t$.  That is, while $G_\eta(t)$ is a nonlinear function of $t$, it is nevertheless the restriction to the boundary of the solution to a system (\ref{model}) of linear equations with respect to independent variables including $t$.  By the same token, $G_\eta(t)$ is a non-linear function of the transformed values $w$ of $\eta$, but it is nevertheless the pushforward of the solution to a system of linear equations in variables $w$ and $x$.   

This suggests that within the broader context of inverse problems for PDEs, it may make sense to ask whether the nonlinear dependence of data on coefficients is governed by some as-yet-unspecified system of linear equations.

\subsection{Resolution of singularities \label{sec-resolution}}

Whereas the piecewise constant structure of layered media has traditionally been analyzed essentially by patching together solutions for individual layers, the global model $\sss_n$ and its dual provide a new framework in which this is no longer necessary.   The proofs of the present paper's main results show that the system $\sss_n$ can be solved directly by separation of variables, leading to an explicit solution to the original model (\ref{model}) in a form suited to time-limited data.  This is possible at least in part because the global model itself has nonsingular coefficients.   

In analogy with the corresponding notion from algebraic geometry, the global model serves as a resolution of singularities for the original model.  Having found one such example, it is natural to ask whether there are others.  More generally, is there a systematic way to associate global models with other classes of PDEs having singular coefficients?

\subsection{The scattering disk\label{sec-mysterious}}

The hybrid laplacian and its eigenfunctions are central to the formulation of both the global model and its solution. The essential problem of computing the Fourier expansion of the backward recurrence formula, as outlined in \S~\ref{sec-essential-questions}, is in the end solved by tensor products of scattering polynomials.  Thus the scattering disk, with its unusual geometry part way between hyperbolic and euclidean, emerges as a key mathematical structure underpinning the relation between data and coefficients for scattering in layered media.  Its role is unexpected and somewhat mysterious, as there seems to be no hint of it in the classical model (\ref{model}) for layered media.


\begin{thebibliography}{10}

\bibitem{BlStSy:2013}
K.~D. Blazek, C.~Stolk, and W.~W. Symes.
\newblock A mathematical framework for inverse wave problems in heterogeneous
  media.
\newblock {\em Inverse Problems}, 29(6):065001, 37, 2013.

\bibitem{BlCoSt:2001}
N.~Bleistein, J.~K. Cohen, and J.~W. Stockwell, Jr.
\newblock {\em Mathematics of multidimensional seismic imaging, migration, and
  inversion}, volume~13 of {\em Interdisciplinary Applied Mathematics}.
\newblock Springer-Verlag, New York, 2001.
\newblock Geophysics and Planetary Sciences.

\bibitem{BrGo:1990}
L.~M. Brekhovskikh and O.~A. Godin.
\newblock {\em Acoustics of Layered Media {I}}, volume~5 of {\em Springer
  Series on Wave Phenomena}.
\newblock Springer, Heidelberg, 1990.

\bibitem{Br:1951}
H.~Bremmer.
\newblock The {W}.{K}.{B}. approximation as the first term of a
  geometric-optical series.
\newblock {\em Comm. Pure Appl. Math.}, 4:105--115, 1951.

\bibitem{BuBu:1983}
K.~P. Bube and R.~Burridge.
\newblock The one-dimensional inverse problem of reflection seismology.
\newblock {\em SIAM Rev.}, 25(4):497--559, 1983.

\bibitem{Bu:1980}
R.~Burridge.
\newblock The {G}elfand-{L}evitan, the {M}archenko, and the {G}opinath-{S}ondhi
  integral equations of inverse scattering theory, regarded in the context of
  inverse impulse-response problems.
\newblock {\em Wave Motion}, 2(4):305--323, 1980.

\bibitem{CoKr:2013}
D.~Colton and R.~Kress.
\newblock {\em Inverse acoustic and electromagnetic scattering theory},
  volume~93 of {\em Applied Mathematical Sciences}.
\newblock Springer, New York, third edition, 2013.

\bibitem{DuXu:2001}
C.~F. Dunkl and Y.~Xu.
\newblock {\em Orthogonal polynomials of several variables}, volume~81 of {\em
  Encyclopedia of Mathematics and its Applications}.
\newblock Cambridge University Press, Cambridge, 2001.

\bibitem{FoGaPaSo:2007}
J.-P. Fouque, J.~Garnier, G.~Papanicolaou, and K.~S{\o}lna.
\newblock {\em Wave propagation and time reversal in randomly layered media},
  volume~56 of {\em Stochastic Modelling and Applied Probability}.
\newblock Springer, New York, 2007.

\bibitem{FuTi:1992}
S.~A. Furman and A.~V. Tikhonravov.
\newblock {\em Basics of Optics of Multilayer Systems}.
\newblock Basics of. Editions Frontires, Gif-sur-Yvette, France, 1992.

\bibitem{Gi:SIAP2014}
P.~C. Gibson.
\newblock The combinatorics of scattering in layered media.
\newblock {\em SIAM J. Appl. Math.}, 74(4):919--938, 2014.

\bibitem{Gi:Dolo2014}
P.~C. Gibson.
\newblock A multivariate interpolation problem arising from the scattering of
  waves in layered media.
\newblock {\em Dolomites Res. Notes Approx. DRNA}, 7:7--15, 2014.

\bibitem{Gr:2008}
L.~Grafakos.
\newblock {\em Classical {F}ourier analysis}, volume 249 of {\em Graduate Texts
  in Mathematics}.
\newblock Springer, New York, second edition, 2008.

\bibitem{Ha:2006}
N.~Hanges.
\newblock {\em Elements of Analytic Hypoellipticity}, volume PM-23 of {\em
  Publica\c{c}\~oes Matem\'aticas}.
\newblock Instituto Nacional de Matem\'atica Pura e Aplicada - IMPA, 2006.

\bibitem{In:2009}
K.~A. Innanen.
\newblock Born series forward modelling of seismic primary and multiple
  reflections: an inverse scattering shortcut.
\newblock {\em Geophysical Journal International}, 177(3):1197--1204, 2009.

\bibitem{Is:1998}
V.~Isakov.
\newblock {\em Inverse problems for partial differential equations}, volume 127
  of {\em Applied Mathematical Sciences}.
\newblock Springer-Verlag, New York, 1998.

\bibitem{Ko:1975}
T.~Koornwinder.
\newblock Two-variable analogues of the classical orthogonal polynomials.
\newblock In {\em Theory and application of special functions ({P}roc.
  {A}dvanced {S}em., {M}ath. {R}es. {C}enter, {U}niv. {W}isconsin, {M}adison,
  {W}is., 1975)}, pages 435--495. Math. Res. Center, Univ. Wisconsin, Publ. No.
  35. Academic Press, New York, 1975.

\bibitem{Ku:1963}
G.~Kunetz.
\newblock Quelques exemples d'analyse d'enregistrements sismiques.
\newblock {\em Geophysical Prospecting}, 11(4):409--422, 1963.

\bibitem{KuFe:2009}
D.~Kurrant and E.~Fear.
\newblock An improved technique to predict the time-of-arrival of a tumor
  response in radar-based breast imaging.
\newblock {\em Biomedical Engineering, IEEE Transactions on}, 56(4):1200
  --1208, april 2009.

\bibitem{Ra:2003}
Rakesh.
\newblock An inverse problem for a layered medium with a point source.
\newblock {\em Inverse Problems}, 19(3):497--506, 2003.

\bibitem{Ra:2008}
Rakesh.
\newblock Inverse problems for the wave equation with a single coincident
  source-receiver pair.
\newblock {\em Inverse Problems}, 24(1):015012, 16, 2008.

\bibitem{Sy:1983}
W.~W. Symes.
\newblock Impedance profile inversion via the first transport equation.
\newblock {\em J. Math. Anal. Appl.}, 94(2):435--453, 1983.

\bibitem{Xu:2015}
Y.~Xu.
\newblock Complex versus real orthogonal polynomials of two variables.
\newblock {\em Integral Transforms Spec. Funct.}, 26(2):134--151, 2015.

\bibitem{Yi:2001}
O.~Yilmaz.
\newblock {\em Seismic data analysis : processing, inversion, and
  interpretation of seismic data}.
\newblock Number~10 in Investigations in Geophysics. Society of Exploration
  Geophysicists, Tulsa, OK, 2nd edition, 2001.
\newblock Edited by Stephen M. Doherty.

\end{thebibliography}
\end{document}